\documentstyle[12pt,amscd]{amsart}
\newtheorem{thm}[equation]{Theorem}
\numberwithin{equation}{section}
\newtheorem{cor}[equation]{Corollary}

\newtheorem{rmk}[equation]{Remark}

\newtheorem{lem}[equation]{Lemma}
\newtheorem{sublem}[equation]{Sublemma}
\newtheorem{conj}[equation]{Conjecture}
\newtheorem{defin}[equation]{Definition}

\newtheorem{prop}[equation]{Proposition}

\begin{document}
\raggedbottom \voffset=-.7truein \hoffset=0truein \vsize=8truein
\hsize=6truein \textheight=8truein \textwidth=6truein
\baselineskip=18truept
\def\mapright#1{\ \smash{\mathop{\longrightarrow}\limits^{#1}}\ }
\def\mapleft#1{\smash{\mathop{\longleftarrow}\limits^{#1}}}
\def\mapup#1{\Big\uparrow\rlap{$\vcenter {\hbox {$#1$}}$}}
\def\mapdown#1{\Big\downarrow\rlap{$\vcenter {\hbox {$\ssize{#1}$}}$}}
\def\mapne#1{\nearrow\rlap{$\vcenter {\hbox {$#1$}}$}}
\def\mapse#1{\searrow\rlap{$\vcenter {\hbox {$\ssize{#1}$}}$}}
\def\mapr#1{\smash{\mathop{\rightarrow}\limits^{#1}}}
\def\ss{\smallskip}
\def\vp{v_1^{-1}\pi}
\def\at{{\widetilde\alpha}}
\def\sm{\wedge}
\def\la{\langle}
\def\ra{\rangle}
\def\on{\operatorname}
\def\spin{\on{Spin}}
\def\kbar{{\overline k}}
\def\qed{\quad\rule{8pt}{8pt}\bigskip}
\def\ssize{\scriptstyle}
\def\a{\alpha}
\def\bz{{\Bbb Z}}
\def\im{\on{im}}
\def\ct{\widetilde{C}}
\def\ext{\on{Ext}}
\def\sq{\on{Sq}}
\def\eps{\epsilon}
\def\ar#1{\stackrel {#1}{\rightarrow}}
\def\br{{\bold R}}
\def\bc{{\bold C}}
\def\bh{{\bold H}}
\def\si{\sigma}
\def\Ebar{{\overline E}}
\def\Sum{\sum}
\def\tfrac{\textstyle\frac}
\def\tb{\textstyle\binom}
\def\Si{\Sigma}
\def\w{\wedge}
\def\equ{\begin{equation}}
\def\b{\beta}
\def\G{\Gamma}
\def\g{\gamma}
\def\psit{\widetilde{\Psi}}
\def\tht{\widetilde{\Theta}}
\def\psiu{{\underline{\Psi}}}
\def\thu{{\underline{\Theta}}}
\def\aee{A_{\text{ee}}}
\def\aeo{A_{\text{eo}}}
\def\aoo{A_{\text{oo}}}
\def\aoe{A_{\text{oe}}}
\def\fbar{{\overline f}}
\def\endeq{\end{equation}}
\def\sn{S^{2n+1}}
\def\zp{\bold Z_p}
\def\A{{\cal A}}
\def\P{{\cal P}}
\def\cj{{\cal J}}
\def\zt{{\bold Z}_2}
\def\bs{{\bold s}}
\def\bof{{\bold f}}
\def\bq{{\bold Q}}
\def\be{{\bold e}}
\def\Hom{\on{Hom}}
\def\ker{\on{ker}}
\def\coker{\on{coker}}
\def\da{\downarrow}
\def\colim{\operatornamewithlimits{colim}}
\def\zphat{\bz_2^\wedge}
\def\io{\iota}
\def\Om{\Omega}
\def\u{{\cal U}}
\def\e{{\cal E}}
\def\exp{\on{exp}}
\def\line{\rule{.6in}{.6pt}}
\def\wbar{{\overline w}}
\def\xbar{{\overline x}}
\def\ybar{{\overline y}}
\def\zbar{{\overline z}}
\def\ebar{{\overline e}}
\def\nbar{{\overline n}}
\def\rbar{{\overline r}}
\def\Ubar{{\overline U}}
\def\et{{\widetilde e}}
\def\ni{\noindent}
\def\coef{\on{coef}}
\def\den{\on{den}}
\def\lcm{\on{l.c.m.}}
\def\vi{v_1^{-1}}
\def\ot{\otimes}
\def\psibar{{\overline\psi}}
\def\mhat{{\hat m}}
\def\exc{\on{exc}}
\def\ms{\medskip}
\def\ehat{{\hat e}}
\def\etao{{\eta_{\text{od}}}}
\def\etae{{\eta_{\text{ev}}}}
\def\dirlim{\operatornamewithlimits{dirlim}}
\def\gt{\widetilde{L}}
\def\lt{\widetilde{\lambda}}
\def\sgd{\on{sgd}}
\def\ord{\on{ord}}
\def\gd{{\on{gd}}}
\def\rk{{{\on{rk}}_2}}
\def\nbar{{\overline{n}}}

\def\N{{\Bbb N}}
\def\Z{{\Bbb Z}}
\def\Q{{\Bbb Q}}
\def\R{{\Bbb R}}
\def\C{{\Bbb C}}
\def\l{\left}
\def\r{\right}
\def\ls{\leq}
\def\gs{\geq}
\def\bg{\bigg}
\def\({\bigg(}
\def\[{\bigg[}
\def\){\bigg)}
\def\]{\bigg]}
\def\colon{{:}\;}
\def\lg{\langle}
\def\rg{\rangle}
\def\t{\text}
\def\f{\frac}
\def\mo{\on{mod}}
\def\vexp{v_1^{-1}\exp}
\def\al{\alpha}
\def\ve{\varepsilon}
\def\bi{\binom}
\def\eq{\equiv}
\def\cs{\cdots}
\def\dstyle{\displaystyle}
\def\Remark{\noindent{\it  Remark}}
\title[Divisibility by 2 and 3 of certain Stirling numbers]
{Divisibility by 2 and 3 of certain Stirling numbers}
\author{Donald M. Davis}
\address{Department of Mathematics, Lehigh University\\Bethlehem, PA 18015, USA}
\email{dmd1@@lehigh.edu}
\date{July 16, 2008}

\keywords{Stirling number, divisibility, James numbers}
\thanks {2000 {\it Mathematics Subject Classification}:
11B73,55Q52.}

\maketitle
\begin{abstract} The numbers $\et_p(k,n)$ defined as $\min(\nu_p(S(k,j)j!):\ j\ge n)$ appear
frequently in algebraic topology.
Here $S(k,j)$ is the Stirling number of the second kind, and $\nu_p(-)$ the exponent of $p$.
Let $s_p(n)=n-1+\nu_p([n/p]!)$.
The author and Sun proved that if $L$ is sufficiently large, then $\et_p((p-1)p^{L}+n-1,n)\ge s_p(n)$.

In this paper, we determine the set of integers $n$ for which $\et_p((p-1)p^{L}+n-1,n)= s_p(n)$ when $p=2$ and when $p=3$.
The condition is roughly that,  in the base-$p$
 expansion of $n$, the sum of two consecutive digits must always be less than $p$.
The result for divisibility of Stirling numbers is, when $p=2$, that for such integers $n$,
$\nu_2(S(2^L+n-1,n))=[(n-1)/2]$.

We also present evidence for conjectures that, if  $n=2^t$ or $2^t+1$, then
the maximum value over all $k\ge n$ of $\et_2(k,n)$ is $s_2(n)+1$.
\end{abstract}

\section{Introduction}\label{intro}
Let $S(k,j)$ denote the Stirling number of the second kind. This satisfies
\begin{equation}\label{Stirdef}S(k,j)j!=(-1)^j\sum_{i=0}^j(-1)^i\tbinom ji i^k.\end{equation}
Let $\nu_p(-)$ denote the exponent of $p$. For $k\ge n$, the numbers $\et_p(k,n)$ defined by
\begin{equation}\label{etdef}\et_p(k,n)=\min(\nu_p(S(k,j)j!):\ j\ge n)\end{equation}
are important in algebraic topology. We will discuss these applications in Section
\ref{topsec}.

In \cite{DS1}, it was proved that, if $L$ is sufficiently large, then
\begin{equation}\label{DS1thm}\et_p((p-1)p^{L}+n-1,n)\ge n-1+\nu_p([n/p]!).\end{equation}
Let $s_p(n)=n-1+\nu([n/p]!)$, as this will appear throughout the paper.
Our main theorems, \ref{mainthm2} and \ref{mainthm3}, give the sets of integers $n$ for which
equality occurs in (\ref{DS1thm}) when $p=2$ and when $p=3$.
Before stating these, we make a slight reformulation to eliminate the annoying $(p-1)p^{L}$.

We define the partial Stirling numbers $a_p(k,j)$ by
$$a_p(k,j)=\sum_{i\not\equiv0\,(p)}(-1)^i\tbinom ji i^k$$
and then
\begin{equation}\label{edef2}e_p(k,n)=\min(\nu_p(a_p(k,j)):\ j\ge n).\end{equation}
Partial Stirling numbers have been studied in \cite{Lun} and \cite{Leng}.

The following elementary and well-known proposition explains the advantage of using $a_p(k,j)$
as a replacement for $S(k,j)j!$: it is that $\nu_p(a_p(k,j))$ is periodic  in $k$. In particular,
$\nu_p(a_p(n-1,n))=\nu_p(a_p((p-1)p^{L}+n-1,n))$ for $L$ sufficiently large, whereas $S(n-1,n)n!=0$.
Thus when using $a_p(-)$, we need not consider the $(p-1)p^{L}$. The second part of the proposition
says that replacing $S(k,j)j!$ by $a_p(k,j)$ merely extends  the numbers $\et_p(k,n)$ for $k\ge n$ in which
we are interested periodically to all integers $k$. An example $(p=3,n=10)$ is given in \cite[p.543]{DSU}.

\begin{prop}\label{elemprop} a. If $t\ge\nu_p(a_p(k,j))$, then
$$\nu_p(a_p(k+(p-1)p^{t},j))=\nu_p(a_p(k,j)).$$
b. If $k\ge n$, then $e_p(k,n)=\et_p(k,n)$.
\end{prop}
\begin{pf} a. (\cite[3.12]{CK}) For all $t$, we have
$$a_p(k+(p-1)p^{t},j)-a_p(k,j)=\sum_{i\not\equiv0\,(p)}(-1)^i\tbinom ji i^k(i^{(p-1)p^t}-1)\equiv0\ (p^{t+1}),$$
from which the conclusion about $p$-exponents is immediate.

b. We have \begin{equation}\label{Sa}(-1)^jS(k,j)j!-a_p(k,j)\equiv 0\ (p^k)\end{equation} since all its terms are multiples of $p^k$. Since
$\et(k,n)\le\nu(S(k,k) k!)<k$, a multiple of $p^k$ cannot affect this value.
\end{pf}

 Our first main result determines the set of values of $n$  for which (\ref{DS1thm}) is sharp when $p=2$.
\begin{thm}\label{mainthm2} For $n\ge1$, $e_2(n-1,n)=s_2(n)$ iff $n=2^{\eps}(2s+1)$ with $0\le\eps\le2$ and $\binom{3s}s$ odd.\end{thm}

\begin{rmk}{\rm Since $\binom{3s}s$ is odd iff binary$(s)$ has no consecutive 1's, another characterization of those $n$
for which $e_2(n-1,n)=s_2(n)$ is those satisfying $n\not\equiv0$ mod 8, and the only consecutive 1's in
binary$(n)$ are, at most, a pair at the end, followed perhaps by one or two 0's.
Alternatively, except at the end, the sum of consecutive bits must be less than 2.}\label{rmk3}\end{rmk}

When $p=3$, the description is  similar.
\begin{thm}\label{mainthm3} Let $T$ denote the set of positive integers for which the sum of two consecutive digits in the
base-$3$ expansion is always less than $3$. Let $T'=\{n\in T:\ n\not\equiv2\,(3)\}$. For integers $a$ and $b$, let
$aT+b=\{an+b:\ n\in T\}$, and similarly for $T'$. Then $e_3(n-1,n)=s_3(n)$ if and only if
$$n\in (3T+1)\cup (3T'+2)\cup (9T+3).$$
\end{thm}

\begin{rmk} {\rm Thus $e_3(n-1,n)=s_3(n)$ iff $n\not\equiv0,6$ (9) and the only consecutive digits in the base-3 expansion of
$n$ whose sum is $\ge3$ are perhaps $\cdots21$, $\cdots12$, or $\cdots210$, each at the very end.}\end{rmk}

The following definition will be used throughout the paper.
\begin{defin} Let $\nbar$ denote the residue of $n$ mod $p$.\end{defin}
\noindent The value of $p$ will be clear from the context. Similarly $\xbar$ denotes the residue of $x$, etc.

\begin{rmk} {\rm As our title suggests, we can interpret our results in terms of divisibility of Stirling numbers.
Suppose $p=2$ or 3
 and $L$ is sufficiently large. The main theorem of \cite{DS1} can be interpreted to say that
\begin{equation}\label{Stirdiv}\nu_p(S((p-1)p^L+n-1,n))\ge (p-1)[\tfrac np]+\nbar-1.\end{equation}
Our main result is that equality occurs in (\ref{Stirdiv}) iff, for $p=3$, $n$ is as in Theorem
 \ref{mainthm3} with $n\not\equiv2$ (9) or, for $p=2$, $n$ is as in Theorem \ref{mainthm2}. We also show that, if $p=3$ and $n=9x+2$, then equality occurs in
 $$\nu_3(S(2\cdot3^L+n-1,n+1))\ge 6x$$
 iff $x\in T'$.}
\end{rmk}

In \cite[(1.5)]{DS2}, a function $T_{k,\a}^p(n,r)$ was introduced, relevant to the proof of (\ref{DS1thm}). We recall it
in Definition \ref{Tdef}.
Useful in our proofs of \ref{mainthm2} and \ref{mainthm3} is the explicit value mod $p$ of $T_{k,2}^p(n,r)$ when $p=2$ and $3$.
(See \ref{Tmod2}, \ref{T12}, and \ref{T1}.)
We obtain this by relating it to $T_{k,1}^p([\frac np], [\frac rp])$ and then evaluating the latter.
This extends \cite[Thm 1.5]{DS2} to the case $\a=1$. Useful in this proof is Theorem \ref{mian}, which is proved in Section \ref{Weissec}
and might be of independent interest.
\begin{defin} If $n$ is a positive integer and $r$ is any integer, let
$$S_1(n,r)=p^{-[\frac{n-1}{p-1}]}\sum_{k\equiv r\, (p)}(-1)^k\tbinom nk,\text{ and }S_2(n,r)=p^{-[\frac{n-1}{p-1}]}\sum_{k\equiv pr\, (p^2)}(-1)^k
\tbinom {pn}k.$$
\end{defin}
\noindent These are integers by \cite{Weis}. They were also studied in \cite{Lun}. The prime $p$ is implicit.
\begin{thm}\label{mian} Let $p$ be an odd prime.
\begin{itemize}
\item[a.] For all $r$, $S_1(n,r)\equiv S_2(n,r)$ mod $p$.
\item[b.] Mod $p$, $S_1(n,r)\equiv\begin{cases}(-1)^{s-1}&\text{if }n=(p-1)s\\
(-1)^{s-1}(\frac {s+1}2+r)&\text{if }n=(p-1)s-1.\end{cases}$
\item[c.] Mod $p$, $S_1(n+p(p-1),r)\equiv -S_1(n,r)$.
\end{itemize}
\end{thm}

Of special interest in algebraic topology is
\begin{equation}\label{ebar}\ebar_p(n):=\max(e_p(k,n):\ k\in\Bbb Z).\end{equation}
In Section \ref{conjsec}, we discuss the relationship between $\ebar_2(n)$,
$e_2(n-1,n)$, and $s_2(n)$. We describe an approach there toward a proof of the following conjecture.
\begin{conj}\label{thm2}If $n=2^t$, then
$$\ebar_2(n)=e_2(n-1,n)=s_2(n)+1,$$
while if $n=2^t+1$, then
$$\ebar_2(n)=e_2(n-1,n)+1=s_2(n)+1.$$
\end{conj}

This conjecture suggests that the inequality
$e_2(n-1,n)\ge s_2(n)$ fails by 1 to be sharp if $n=2^t$, while if $n=2^t+1$, it is sharp but the maximum value of $e_2(k,n)$
occurs for a value of $k\ne n-1$.

\section{Proof of Theorem \ref{mainthm2}}\label{pfthm1sec}
In this section, we prove Theorem \ref{mainthm2}, utilizing results of \cite{DS2} and some work with binomial coefficients.
The starting point is the following result of \cite{DS2}. In this section, we abbreviate $\nu_2(-)$ as $\nu(-)$.
\begin{thm}\label{DS2thm} $(\cite[1.2]{DS2})$ For all $n\ge0$ and $k\ge0$,
$$\nu\left( 2^k k!\sum_i\tbinom n{4i+2}\tbinom ik\right)\ge \nu([n/2]!).$$
\end{thm}
The bulk of the work is in proving the following refinement. The inequality is immediate from \ref{DS2thm}.
\begin{thm} \label{helpthm}Let $n$ be as in Theorem \ref{mainthm2}, and, if $n>4$, define $n_0$ by $n=2^e+n_0$ with $0<n_0<2^{e-1}$. Then
\begin{equation}\label{DSineq} \nu\left(\tbinom {n-1}k 2^k k!\sum_i\tbinom n{4i+2}\tbinom ik\right)\ge \nu([n/2]!)
\end{equation}
for all $k$, with equality if and only if
\begin{equation}\label{3cases}k=\begin{cases}0&1\le n\le4\\ n_0-1&n\not\equiv0\ (\mo\, 4),\ n>4\\
n_0-2&n\equiv0\ (\mo\, 4),\  n>4.\end{cases}\end{equation}
\end{thm}

\begin{pf*}{Proof that Theorem \ref{helpthm} implies the ``if" part of Theorem \ref{mainthm2}}
By (\ref{DS1thm}),  $e_2(n-1,n)\ge s_2(n)$ for all $n$.
Thus it will suffice to prove that if $n$ is as in Theorem \ref{helpthm}, then
\begin{equation}\label{nuaeq}\nu(a_2(n-1,n))=s_2(n).\end{equation}

Note that
$$0=(-1)^nS(n-1,n)n!=-a_2(n-1,n)+\sum\tbinom n{2k}(2k)^{n-1}.$$
Factoring $2^{n-1}$ out of the sum shows that
(\ref{nuaeq}) will follow from showing
\begin{equation}\label{noL} \sum \tbinom n{2k}k^{n-1}=\nu([n/2]!).\end{equation}

The sum in (\ref{noL}) may be restricted to odd values of $k$, since terms with even $k$ are more 2-divisible than
the claimed value. Write $k=2j+1$ and apply the Binomial Theorem, obtaining
\begin{equation}\label{Stirl}\sum_j\tbinom n{4j+2}\sum_\ell 2^\ell j^\ell\tbinom{n-1}\ell=\sum_j\tbinom n{4j+2}\sum_\ell
2^\ell\tbinom{n-1}\ell\sum_iS(\ell,i)i!\tbinom ji.\end{equation}
Here we have used the standard fact that $j^{\ell}=\sum S(\ell,i)i!\binom ji$.

Recall that $S(\ell,i)=0$ if $\ell<i$, and $S(i,i)=1$. Terms in the right hand side of (\ref{Stirl}) with $\ell=i$ yield
$$\sum_i\tbinom{n-1}i 2^ii!\sum_j\tbinom n{4j+2}\tbinom ji,$$
which we shall call $A_n$.
By Theorem \ref{helpthm}, if $n$ is as in Theorem \ref{mainthm2}, $\nu(A_n)=\nu([n/2]!)$ since all $i$-summands have
2-exponent $\ge\nu([n/2]!)$, and exactly one of them has 2-exponent equal to $\nu([n/2]!)$. Terms in (\ref{Stirl})
with $\ell>i$ satisfy $$\nu(\text{term})>\nu\left(2^ii!\sum_j\tbinom n{4j+2}\tbinom ji\right),$$ the RHS of which
is $\ge\nu([n/2]!)$ by \ref{DS2thm}. The claim (\ref{noL}), and hence Theorem \ref{mainthm2}, follows.
\end{pf*}

We recall the following definition from \cite[1.5]{DS2}.
\begin{defin}\label{Tdef} Let $p$ be any prime.
 For $n,\a,k\ge0$ and $r\in\Bbb Z$, let
$$T^p_{k,\a}(n,r):=\frac{k!p^k}{[n/p^{\a-1}]!}\sum_i(-1)^{p^\a i+r}\binom n{p^\a i+r}\binom ik.$$
\end{defin}
\noindent In the remainder of this section, we have $p=2$ and omit writing it as a superscript of $T$.

By \ref{DS2thm}, Theorem \ref{helpthm} is equivalent to the following result, to the proof of which the rest of this section will be devoted.
\begin{thm}\label{newthm} If $n$ is as in Theorem \ref{helpthm}, then $\binom{n-1}kT_{k,2}(n,2)$ is odd if and only if $k$ is as in $(\ref{3cases})$.\end{thm}

Central to the proof of \ref{newthm} is the following result, which will be proved at the end of this section.
This result applies to all values of $n$, not just those as in Theorem \ref{helpthm}. This result is the complete
evaluation of $T_{k,2}(n,2)$ mod $2$.
\begin{thm}\label{Tmod2} If $4k+2>n$, then $T_{k,2}(n,2)=0$. If $4k+2\le n$, then, mod $2$,
$$T_{k,2}(n,2)\equiv\binom{[n/2]-k-1}{[n/4]}.$$
\end{thm}

\begin{pf*}{Proof of Theorem \ref{newthm}}
The cases $n\le 4$  are easily verified and not considered further.

First we establish that $\binom{n-1}kT_{k,2}(n,2)$ is odd for the stated values of $k$. We have
$$\tbinom{n-1}k=\begin{cases}\binom{2^e+n_0-1}{n_0-1}&\text{if }n_0\not\equiv0\ (\mo \, 4)\\
\binom{2^e+n_0-1}{n_0-2}&\text{if }n_0\equiv0\ (\mo\, 4),\end{cases}$$
which is clearly odd in both cases. Here and throughout we use the well-known fact that, if $0\le\eps_i,\delta_i\le p-1$, then
\begin{equation}\label{Lucas}\binom{\sum \eps_ip^i}{\sum \delta_ip^i}\equiv\prod \binom{\eps_i}{\delta_i}\ (\mo\, p).\end{equation}

Now we show that $T_{k,2}(n,2)$ is odd when $n$ and $k$ are as \ref{mainthm2} and (\ref{3cases}).

\noindent {\bf Case 1}: $n_0=8t+4$ with $\binom{3t}t$ odd, and $k=8t+2$. Using \ref{Tmod2}, with all equivalences mod 2,
$$T_{k,2}(n,2)\equiv\binom{2^{e-1}+4t+2-(8t+2)-1}{2^{e-2}+2t+1}\equiv\binom{-4t-1}{2t+1}\equiv\binom{6t+1}{2t+1}\equiv\binom{3t}t.$$
{\bf Case 2}: $n_0=4t+\eps$, $\eps\in\{1,2\}$, $\binom{3t}t$ odd, $k=4t+\eps-1$. Then
$$T_{k,2}(n,2)\equiv\binom{2^{e-1}+2t+\eps-1-(4t+\eps-1)-1}{2^{e-2}+t}\equiv\binom{-2t-1}t\equiv\binom{3t}t.$$
{\bf Case 3}: $n_0=4t+3$, $\binom{3(2t+1)}{2t+1}$ odd, $k=4t+2$. Then
$$T_{k,2}(n,2)\equiv\binom{2^{e-1}+2t+1-(4t+2)-1}{2^{e-2}+t}\equiv\binom{-2t-2}t\equiv\binom{3t+1}t\equiv\binom{2(3t+1)+1}{2t+1}.$$

Now we must show that, if $n$ is as in Theorem \ref{mainthm2} and $k$ does not have the value specified in
(\ref{3cases}), then $\binom{n-1}kT_{k,2}(n,2)$ is even. The notation of Theorem \ref{helpthm} is
continued. We divide into cases.

{\bf Case 1:} $k\ge n_0$. Here $\binom{n-1}k$ odd implies $k\ge 2^e$, but then $4k+2>n$ and so by Theorem \ref{Tmod2}, $T_{k,2}(n,2)=0$.
 Hence $\binom{n-1}k T_{k,2}(n,2)$ is even.

\ss
{\bf Case 2:} $n_0=4t+4$, $k=n_0-1$. Here $T_{k,2}(n,2)\equiv\binom{-(2t+2)}{t+1}\equiv\binom{3t+2}{t+1}$. If $t$ is even, this is even, and if
$t=2s-1$, this is congruent to $\binom{3s-1}s$ which is even, since if $\nu(s)=w$, then $2^w\not\in 3s-1$; i.e.,
the decomposition of $3s-1$ as a sum of distinct 2-powers does not contain $2^w$.
\ss

{\bf Case 3:} $n_0=4t+\eps$, $1\le \eps\le3$, and $k<n_0-1$. Here $$\binom{n-1}kT_{k,2}(n,2)\equiv\binom{4t+\eps-1}k
\binom{2^{e-1}+2t+[\eps/2]-k-1}{2^{e-2}+t}.$$
If $k\le 2t+[\eps/2]-1$, then the second factor is even due to the $i=e-2$ factor in (\ref{Lucas}). If $k>2t+[\eps/2]-1$, the second factor
is congruent to $\binom{-(k+1-2t-[\eps/2])}t\equiv\binom{k-t-[\eps/2]}t$. For
$\binom{4t+\eps-1}k\binom{k-t-[\eps/2]}t$ to be odd would require one of the following:
\begin{eqnarray*} &&\eps=1,\ k=4i,\text{ and }\tbinom ti\tbinom{4i-t}t\text{ odd}\\
&&\eps=2,\ k=4i+\langle0,1\rangle,\text{ and }\tbinom ti\tbinom{4i-t-\langle 1,0\rangle}t\text{ odd.}\\
&&\eps=3,\ k=4i+\langle0,2\rangle,\ \tbinom ti\tbinom{4i-t+\langle-1,1\rangle}t\text{ odd.}
\end{eqnarray*}
But all these products are even if $i<t$ by Lemma \ref{comblem}. If $i=t$, since $k<n_0-1$, we obtain a $\binom{3t-1}t$ factor, which is even, as in Case 2.

\ss
{\bf Case 4:} $n_0=4t+4$ and $k<n_0-2$. Note that $t$ must be even since $n\not\equiv0\ (8)$ in \ref{helpthm}. We have
$$\binom{n-1}kT_{k,2}(n,2)\equiv\binom{4t+3}k\binom{2^{e-1}+2t+1-k}{2^{e-2}+t+1}.$$
The case $k\le2t+1$ is handled as in Case 3. If $k> 2t+1$, then, similarly to Case 3, it reduces to
$\binom{4t+3}k\binom{k-t-1}{t+1}$. If $k=4t$ or $4t+1$, then we obtain $\binom{3t-1}{t+1}$ or $\binom{3t}{t+1}$, which are even since $t$ is even.
Now suppose $k=4i+\Delta$ with $0\le\Delta\le3$ and $i<t$. Since $t$ is even, if $\Delta$ is odd, then $\binom{k-t-1}{t+1}$ is even. For $\Delta=0$ or 2, we obtain $\binom ti\binom{4i-t\pm1}{t+1}$. Since $t$ is even,
we use $\binom{2A+1}{2B+1}\equiv\binom{2A}{2B}$ to obtain $\binom ti\binom{4i-t-\langle0,2\rangle}t$, which is even by Lemma \ref{comblem}.
\end{pf*}

The following result implies the ``only if" part of Theorem \ref{mainthm2}.
\begin{thm} Assume $n\equiv0$ mod $8$ or $n=2^{\eps}(2s+1)$ with $0\le\eps\le2$ and $\binom{3s}s$ even. Then for all $N\ge n$, we have
$\nu_2(a_2(n-1,N))> s_2(n)$.\end{thm}
\begin{pf}
Combining aspects of \ref{helpthm}, \ref{Tmod2}, and \ref{phigen}, the theorem will follow from showing that for $n$ as in the theorem
and $N\ge n$ satisfying $[N/4]=[n/4]$, we have
\begin{equation}\label{nN} \sum_{4k+2\le N}\binom{n-1}k\binom{[N/2]-k-1}{[N/4]}\equiv0\ (2).\end{equation}
Note that if $[N/4]>[n/4]$, then $\frac{[N/2]!}{[n/2]!}$ is even in the 2-primary analogue of the proof of \ref{phigen}.

When $n=8\ell$, it is required to show that $\sum\binom{8\ell-1}k\binom{4\ell-k-1}{2\ell}$ and $\sum\binom{8\ell-1}k\binom{4\ell-k}{2\ell}$
are both even. The first corresponds to $N=n$ or $n+1$, and the second to $N=n+2$ or $n+3$. The first is proved by noting easily that the
summands for $k=2j$ and $2j+1$ are equal. The second follows from showing that the summands for $k=2j$ and $2j-1$ are equal. This is easy unless
$2j=8i$. For this, we need to know that $\binom{2\ell-4i}{\ell}\binom\ell i$ is always even, and this follows easily from showing that
the binary expansions of $\ell-4i$, $\ell-i$, and $i$ cannot be disjoint.

For $n=2^\eps(2s+1)$ with $\binom{3s}s$ even, all summands in (\ref{nN})
can be shown to be even when $n=2^e+n_0$ with $0<n_0<2^{e-1}$ and $N=n$ using the proof of
Theorem \ref{newthm}. For such $n$ and $N>n$, the main case to consider is $n=8a+4$ and $N=n+2$. Then we need
$\binom{8a+3}k\binom{4a+2-k}{2a+1}\equiv0$ mod 2. For this to be false, $k$ must be odd. But then we have
$$\tbinom{8a+3}k\tbinom{4a+2-k}{2a+1}\equiv\tbinom{8a+3}{k-1}\tbinom{4a+1-(k-1)}{2a+1}\equiv0$$
by the result for $N=n$ with $k$ replaced by $k-1$.

If $n=2^{e+d}+\cdots+2^e+n_0$ with $d>0$ and $0<n_0<2^{e-1}$, then (\ref{nN}) for $n=N$ is proved when $k$ does not have the special value of
(\ref{3cases}) just as in the second part of the proof of  \ref{newthm}. We illustrate what happens when $k$ does have the special value
by considering what happens to Case 1 just after (\ref{Lucas}). The binomial coefficient there becomes
$$\binom{2^{e+d-1}+\cdots+2^{e-1}-4t-1}{2^{e+d-2}+\cdots+2^{e-2}+2t+1},$$
which is 0 mod 2 by consideration of the $2^{e-1}$ position in (\ref{Lucas}). For $N>n$, the argument is essentially the same as that of the previous
paragraph.
\end{pf}

The following lemma was used above.
\begin{lem}\label{comblem} Let $i<t$, $-2\le\delta\le1$, and if $\delta=-2$, assume that $t$ is even. Then $\binom ti\binom{4i-t+\delta}t$ is even.
\end{lem}
\begin{pf} Assume that $\binom ti\binom{4i-t+\delta}t$ is odd. Then $i$, $t-i$, and $4i-2t+\delta$ have disjoint binary expansions. If $\delta=0$ or 1, then letting $\ell=t-i$ and $r=2i-t$, we infer that $\ell+r$, $\ell$, and $2r$ are disjoint with $\ell$ and $r$ positive, which is impossible by Sublemma \ref{subl}.2. If $\delta=-1$ and $t$ is odd, then two of $i$, $t-i$, and $4i-2t-1$ are odd, and so they cannot be disjoint. Thus we may assume $t$ is even and $\delta=-1$ or $-2$. Let $\ell=t-i$ and $r=2i-t-1$. Then
$\ell+r+1$, $\ell$, and $2r$ are disjoint with $\ell$ and $r$ positive and $r$ odd, which is impossible by
Sublemma \ref{subl}.3.
\end{pf}
\begin{sublem}\label{subl} Let $\ell$ and $r$ be nonnegative integers.
\begin{enumerate} \item Then $\ell$, $2r+1$, and $\ell+r+1$ do not have disjoint binary expansions.
\item If $\ell$ and $r$ are positive, then $\ell$, $2r$, and $\ell+r$ do not have disjoint binary expansions.
\item If $\ell$ is positive and $r$ is odd, then $\ell$, $2r$, and $\ell+r+1$ do not have disjoint binary expansions.
\end{enumerate}
\end{sublem}
\begin{pf} (1) Assume that $\ell$ and $r$ constitute a minimal counterexample. We must have $\ell=2\ell'$ and $r=2r'+1$.
Then $\ell'$ and $r'$ yield a smaller counterexample.

(2) Assume that $\ell$ and $r$ constitute a minimal counterexample. If $r$ is even, then $\ell$ must be even, and so dividing each by 2 gives a smaller counterexample. If $r=1$, then $\ell$, 2, and $\ell+1$ are disjoint, which is
impossible, since the only way for $\ell$ and $\ell+1$ to be disjoint is if $\ell=2^e-1$. If $r=2r'+1$ with $r'>0$,
and $\ell=2\ell'$, then $\ell'$ and $r'$ form a smaller counterexample. If $r=2r'+1$ and $\ell=2\ell'+1$,
then $\ell'$, $2r'+1$, and $\ell'+r'+1$ are disjoint, contradicting (1).

(3) Let $r=2r'+1$. Then $\ell$ must be even ($=2\ell'$). Then $\ell'$, $2r'+1$, and $\ell'+r'+1$ are disjoint,
contradicting (1).
\end{pf}

The following lemma together with Theorem \ref{T12} implies Theorem \ref{Tmod2}. Its proof uses the following definition,
which will be employed throughout the paper.
\begin{defin} \label{ddef} Let $d_p(-)$ denote the number of 1's in the $p$-ary expansion.\end{defin}
\begin{lem}\label{T1lem} Mod 2,
\begin{equation}\label{T1eq}T_{k,1}(n,r)\equiv\begin{cases}\binom{n-k-1}{[(n-1+\rbar)/2]}&n>k\\
0&n\le k\end{cases}\end{equation}
\end{lem}
\begin{pf} The proof is by induction on $k$. Let $f_k(n,r)$ denote the RHS of (\ref{T1eq}) mod 2.
It is easy to check that $f_0(n,r)=\delta_{d_2(n),1}$, agreeing with $T_{0,1}(n,r)$ as determined in
(\ref{T1del}). Here and throughout $\delta_{i,j}$ is the Kronecker function. From Definition \ref{Tdef}, mod 2,
$T_{k,1}(1,r)\equiv\delta_{k,0}$.
This is what causes the dichotomy in (\ref{T1eq}).

 By \cite[(2.3)]{DS2},  if $k>0$, then
\begin{equation}\label{indnform}T_{k,1}(n,r)+rT_{k-1,1}(n,r+2)=-T_{k-1,1}(n-1,r+1).\end{equation}
Noting that $f$ only depends on the mod 2 value of $r$, the lemma
follows from
\begin{eqnarray*}f_k(n,0)&=&f_{k-1}(n-1,1)\\
f_k(n,1)&=&f_{k-1}(n,1)+f_{k-1}(n-1,0),\end{eqnarray*}
which are immediate from the definition of $f$ and Pascal's formula.
\end{pf}

\section{Mod $p$ values of $T$-function}\label{Weissec}
We saw in Theorem \ref{newthm} that knowledge of the mod 2 value of the $T$-function of \cite{DS2}
played an essential role in proving Theorem \ref{mainthm2}. A similar situation occurs when $p=3$.
The principal goal of this section is the determination of $T_{k,2}^3(n,r)$, obtained by combining
Theorems \ref{T12} and \ref{T1}. We also prove Theorem \ref{mian}, which is used in the proof of \ref{T12},
but may be of intrinsic interest.

We begin by recording a well-known proposition which will be used throughout the paper.

\begin{prop}\label{Lagra} If $n\ge0$, then $\nu_p(n!)=\frac1{p-1}(n-d_p(n))$, and hence $\nu_p(\binom nb)=\frac1{p-1}(d_p(b)+d_p(n-b)-d_p(n))$.
\end{prop}

The following result extends \cite[Thm 1.5]{DS2} to include the case $\a=1$.
\begin{thm}\label{T12} Let $p$ be any prime. For any $\a\ge1$, we have the congruence, mod $p$,
$$T^p_{k,\a+1}(n,r)\equiv(-1)^{\rbar}\binom{\nbar}{\rbar}T^p_{k,\a}([\tfrac np],[\tfrac rp]).$$
\end{thm}
\begin{pf} This was proved for $\a\ge2$ in \cite[Thm 1.5]{DS2}. The only place that the proof
of that result does not work when $\a=1$ is in the initial step of Case 3 of \cite[p.5548]{DS2}. Required to complete that proof is
$$T^p_{0,2}(pn,pr)\equiv T^p_{0,1}(n,r) \pmod p.$$
This just says, mod $p$,
\begin{equation}\label{Kron}{\tfrac 1{n!}}\sum_{i\equiv pr\,(p^2)}(-1)^i\tbinom {pn}i\equiv{\tfrac 1{n!}}\sum_{i\equiv r\,(p)}(-1)^i\tbinom ni.
\end{equation}
When $p$ is odd, this follows immediately from part a of Theorem \ref{mian}, since $\nu_p(n!)\le[(n-1)/(p-1)]$.

We prove (\ref{Kron}) when $p=2$ by showing that both sides equal $\delta_{1,d_2(n_0)}$. The RHS equals
\begin{equation}\label{T1del}\frac1{2^{n-d_2(n)}u}\cdot2^{n-1}
\equiv 2^{d_2(n)-1}\equiv\delta_{d_2(n),1}\,(\mo \ 2),\end{equation}
with $u$ odd, while the LHS is
$\dstyle{{\tfrac 1{n!}}\sum_{i\equiv2r\, (4)}\tbinom{2n}i\equiv\frac1{2^{n-d_2(n)}}\cdot\begin{cases}
2^{2n-2}&n\text{ odd}\\2^{n-1}u'&n\text{ even,}\end{cases}}$,
and this also equals $\delta_{d_2(n),1}$. Here we have used
$$\sum_{i\equiv r\, (4)}\tbinom ni=2^{n-2}+\eps_{n,r}2^{[n/2]-1},\text{ with }
\eps_{n,r}=\begin{cases}0&n-2r\equiv 2\ (\mo\,4)\\ 1&n-2r\equiv -1,0,1\ (\mo\,8)\\ -1&n-2r\equiv3,4,5\ (\mo\,8),\end{cases}$$
which is easily proved by induction on $n$.
\end{pf}

Next we discuss Theorem \ref{mian} and give its proof.
 First we note that the definitions of $S_1$ and $S_2$ in it are similar to \cite[(3.4)]{DS2}, but differ regarding the role of the
second variable in $S_2$.
 We remark that part b of \ref{mian} was given by Lundell in \cite{Lun}, although he merely said ``the proof is a straightforward but
somewhat tedious induction." Part a is  of particular interest to us.

\begin{pf*}{Proof of Theorem \ref{mian}} Throughout this proof, $p$ denotes an odd prime.
We will work with polynomials in the ring $R:={\Bbb F}_p[x]/(x^p-1)$. In $R$, let
\begin{equation}\label{pqdef}P_n(x)=\sum_{r=0}^{p-1}S_1(n,r)x^r\text{ and }Q_n(x)=\sum_{r=0}^{p-1}S_2(n,r)x^r.\end{equation}
Also in $R$, let
$$\psi(x)=\frac{(1-x)^p-(1-x^p)}{p(1-x)}=\frac{(1-x)^{p-1}-(1+\cdots+x^{p-1})}p.$$

We will prove later the following result, which immediately implies part a.
\begin{thm} \label{PQthm}For $1\le d\le p-1$ and $m\ge0$, we have in $R$
$$P_{(p-1)m+d}(x)=\psi(x)^m(1-x)^d=Q_{(p-1)m+d}(x).$$
\end{thm}
Parts b and c follow from Theorem \ref{PQthm} and the following result, which we will also
prove later. The numbering of the parts is related to the corresponding part of Theorem \ref{mian}.
\begin{lem}\label{philem} We have, in $R$,
\begin{itemize}
\item[b.i.] $(\psi(x)+1)(1-x)^{p-1} = 0$,
\item[b.ii.] $(\psi(x)+x^{(p-1)/2})(1-x)^{p-2}=0$, and
\item[c.] $\psi(x)^p=-1$.
\end{itemize}
\end{lem}
The deduction of \ref{mian}.bc is straightforward. For the first part of b, we have in $R$
\begin{eqnarray*}&&(-1)^{s-1}(1+x+\cdots+x^{p-1})\\
&=&(-1)^{s-1}(1-x)^{p-1}\\
&=&\psi(x)^{s-1}(1-x)^{p-1}\\
&=&P_{(p-1)s}(x)\\
&=&\sum_{r=0}^{p-1}S_1((p-1)s,r)x^r.
\end{eqnarray*}

Noting that $(1-x)^{p-2}=1+2x+3x^2+\cdots+(p-1)x^{p-2}$, the second part of \ref{PQthm}.b follows from
the following analysis of coefficients of polynomials in $R$.
\begin{eqnarray*}&&S_1((p-1)(s-1)+p-2,r)\\
&=&[x^r]P_{(p-1)(s-1)+p-2}(x)\\
&=&[x^r](\psi(x)^{s-1}(1-x)^{p-2})\\
&=&[x^r]((-1)^{s-1}x^{(s-1)(p-1)/2}(1-x)^{p-2})\\
&=&(-1)^{s-1}[x^{r+(s-1)/2}](1-x)^{p-2}\\
&=&(-1)^{s-1}(r+(s-1)/2+1).\end{eqnarray*}
Note that exponents of $x$ may be considered mod $p$.
The deduction of \ref{mian}.c from \ref{philem}.c is much easier, and omitted.
\end{pf*}
\begin{pf*}{Proof of Theorem \ref{PQthm}} We first show the theorem is true when $m=0$. The argument for $P$
is similar to, and easier than, the following argument for $Q$. Let $1\le d\le p-1$. Note that, mod $p$,
$$Q_d(x^p)=\sum_r S_2(d,r)x^{pr}=\sum_r(-1)^{pr}\tbinom{pd}{pr}x^{pr}\equiv\sum_{r=0}^d(-1)^r\tbinom d rx^{pr}=(1-x^p)^d.$$
Thus the same is true when $x^p$ is replaced by $x$. Note that
here we are dealing with polynomials mod $p$, but not in the ring $R$ used earlier.

Next we prove that for any $n\not\equiv1$ mod $(p-1)$
\begin{equation}\label{Pphi}P_{n+p-1}(x)=\psi(x)P_n(x)\end{equation}
in $R$. To see this, first note that if $n\not\equiv1$ mod $(p-1)$,
\begin{equation}\label{S1}S_1(n,r)=S_1(n-1,r)-S_1(n-1,r-1).\end{equation}
Note that the need for $n\not\equiv1$ is so that $[(n-1)/(p-1)]=[(n-2)/(p-1)]$.
Similarly, for $n\not\equiv1$ mod $(p-1)$
\begin{equation}\label{S1p-1}S_1(n+p-1,r)={\tfrac 1p}\sum_{i=0}^p(-1)^i\tbinom pi S_1(n-1,r-i).\end{equation}
Since $S_1(n-1,r)=S_1(n-1,r-p)$, this becomes
\begin{eqnarray} S_1(n+p-1,r)&=&\sum_{i=1}^{p-1}(-1)^i\tfrac1p\tbinom piS_1(n-1,r-i)\nonumber\\
&=&\sum_{i=1}^{p-2}\a_iS_1(n,r-i),\label{p-2}\end{eqnarray}
where
\begin{equation}\label{phialph}\psi(x)=\sum_{i=1}^{p-2}\a_ix^i.\end{equation} At the last step, we have used (\ref{S1}). The equation (\ref{p-2})
translates to (\ref{Pphi}).

A similar argument, sketched below, shows that for any $n\not\equiv1$ mod $(p-1)$
\begin{equation}\label{Qphi}Q_{n+p-1}(x)=\psi(x)Q_n(x)\end{equation}
in $R$. The $S_2$-analogue of (\ref{S1}) is true mod $p$, obtained from
$$S_2(n,r)=S_2(n-1,r)+p^{-[\frac{n-1}{p-1}]}\sum_{i=1}^{p-1}(-1)^i\tbinom pi\sum_{k\equiv pr-i\ (p^2)}
(-1)^k\tbinom{p(n-1)}k -S_2(n-1,r-1)$$
by noting that the $k$-sums are divisible by $p^{[(n-1)/(p-1)]}$ by \cite{Weis}, and so since
$\binom pi\equiv0$ mod $p$, then each $i$-summand is 0 mod $p$. The $S_2$-analogue of (\ref{S1p-1}), mod $p$, is obtained
similarly, using that $(1-x)^{p^2}\equiv (1-x^p)^p$ mod $p$. The argument for (\ref{Qphi}) is completed
as in (\ref{p-2}).

Theorem \ref{PQthm} with $d\ne1$ is immediate from (\ref{Pphi}) and (\ref{Qphi}) plus the validity when $m=0$
established in the first paragraph of this proof. The proof when $d=1$ requires the following three lemmas.
\begin{lem}\label{lem1} If $n$ is odd and $n-2r\equiv0$ mod $p$, then $S_1(n,r)=0=S_2(n,r)$.\end{lem}
\begin{pf} Since $n-r\equiv r$ mod $p$, both $\binom nr$ and $\binom n{n-r}$ occur in the sum for $S_1(n,r)$, and with opposite sign
since $n$ is odd. Hence all terms in the sum occur in cancelling pairs. The same is true of all terms in the sum for $S_2(n,r)$
since $pn-pr\equiv pr$ mod $p^2$.\end{pf}
\begin{lem}\label{lem2} If $(1-x)f(x)=0$ in $R$, then $f(x)=c(1+x+\cdots+x^{p-1})$ for some $c$.\end{lem}
\begin{pf} Let $f(x)=c_0+c_1x+\cdots+c_{p-1}x^{p-1}$. The given equation implies $c_0=c_1=\cdots=c_{p-1}$.\end{pf}
\begin{lem}\label{lem3} For $t\in\Z$, let $R_t\subset R$ denote the span of $x^i-x^{t-i}$ for all $i$. If $g(x)\in R_t$, then
$g(x)\psi(x)\in R_{t-1}$.\end{lem}
\begin{pf} Since $\psi(x)$ is a linear combination of various $x^j+x^{p-1-j}$, the lemma follows from the observation that
$$(x^i-x^{t-i})(x^j+x^{-1-j})=x^{i+j}-x^{t-1-i-j}+x^{i-1-j}-x^{t+j-i}.$$
\end{pf}
Note that if $g(x)\in R_t$, then $[x^{t/2}]g(x)=0$.

Now we prove the case $d=1$ of Theorem \ref{PQthm}. We have
$$P_{(p-1)m+1}(x)\cdot(1-x)=P_{(p-1)m+2}(x)=(1-x)^2\psi(x)^m.$$
By Lemma \ref{lem2}, $\Delta_m(x):=P_{(p-1)m+1}(x)-(1-x)\psi(x)^m$ has all coefficients equal.
By Lemma \ref{lem1}, if $(p-1)m+1-2r\equiv0$ mod $p$, then $[x^r]P_{(p-1)m+1}(x)=0$. Note that here $r=(1-m)/2$, with exponents always
considered mod $p$ in $R$.
In the notation of Lemma \ref{lem3}, $1-x\in R_1$, and hence by that lemma, $(1-x)\psi(x)^m\in R_{1-m}$. Thus $[x^{(1-m)/2}]((1-x)\psi(x)^m)=0$.
Thus $[x^{(1-m)/2}]\Delta_m(x)=0$, and hence $\Delta_m(x)=0$, as desired.
\end{pf*}

\begin{pf*}{Proof of Lemma \ref{philem}} To prove b.i., we prove $\psi(x)+1$ is divisible by $(1-x)$ by showing $\psi(1)\equiv-1$ mod $p$.
Note that $\sum_{i=0}^{p-1}((-1)^i\binom{p-1}i-1)=-p$, and hence
$$\psi(1)=\tfrac 1p\sum_{i=0}^{p-1}((-1)^i\tbinom{p-1}i-1)=-1.$$

 To prove b.ii., we prove $g(x):=\psi(x)+x^{(p-1)/2}$ is divisible by $(1-x)^2$. Since $g(1)=0$, it remains to show that the derivative
 satisfies $g'(1)=0$; i.e., that $\psi'(1)+\frac{p-1}2\equiv0$ mod $p$. Let $\a_i=\frac 1p((-1)^i\binom{p-1}i-1)$. Then $\psi'(1)=\sum_{i=1}^{p-1}i\a_i$.
 Since
 $$-(p-1)(1-x)^{p-2}=\tfrac d{dx}(1-x)^{p-1}=\sum_{i=1}^{p-1}(-1)^i\tbinom{p-1}iix^{i-1},$$
 setting $x=1$ shows $\sum_{i=1}^{p-1}(-1)^i\tbinom{p-1}ii=0$ and thus
 $$p\psi'(1)=\sum_{i=1}^{p-1}pi\a_i=\sum_{i=1}^{p-1}((-1)^i\tbinom{p-1}i-1)i=-\sum_{i=1}^{p-1}i=-\tfrac{p(p-1)}2,$$
 and hence $\psi'(1)+\frac{p-1}2=0$, as desired.

 To prove c, we use $x^p=1$, $(A+B)^p\equiv A^p+B^p$, and $i^p\equiv i$, and obtain, in $R$,
 $$\psi(x)^p=\sum_{i=0}^{p-1}\frac{(-1)^i\tbinom{p-1}i-1}p=0-1.$$
 \end{pf*}

Now we give the mod 3 values of $T_{k,1}^3(-,-)$. The mod 3 values of $T_{k,2}^3(-,-)$ can be obtained from this using Theorem \ref{T12}.
Throughout the rest of this section and the next, the superscript 3 on $T$ is implicit.
\begin{thm}\label{T1} Let $n=3m+\delta$ with $0\le\delta\le2$.
\begin{itemize}
\item If $n-k=2\ell$, then, mod $3$, $T_{k,1}(n,r)$ is given by

\begin{center}
\begin{tabular}{cc|ccc}
&&&$\delta$&\\
&&$0$&$1$&$2$\\
\hline
&$0$&$\tbinom{\ell-1}{m-1}$&$\tbinom{\ell-1}m$&$-\tbinom{\ell-1}m$\\
$r\pmod 3$&&&&\\
&$1,2$&$-\tbinom{\ell-1}m$&$\tbinom{\ell-1}m$&$-\tbinom{\ell-1}m$
\end{tabular}
\end{center}
\item If $n-k=2\ell+1$, then, mod $3$, $T_{k,1}(n,r)$ is given by

\begin{center}
\begin{tabular}{cc|ccc}
&&&$\delta$&\\
&&$0$&$1$&$2$\\
\hline
&$0$&$0$&$\tbinom\ell m$&$0$\\[3pt]
$r \pmod  3$&$1$&$\tbinom\ell m$&$-\tbinom\ell m$&$0$\\[3pt]
&$2$&$-\tbinom\ell m$&$0$&$0$
\end{tabular}
\end{center}
\end{itemize}
\end{thm}
\begin{pf} By \cite[(2.3)]{DS2}, we have \begin{equation}\label{recur}T_{k,1}(n,r)+rT_{k-1,1}(n,r+3)=-T_{k-1,1}(n-1,r+2),\end{equation}
yielding an inductive determination of $T_{k,1}$ starting with $T_{0,1}$.
One can verify that the mod 3 formulas of Theorem \ref{T1} also satisfy
(\ref{recur}). For example, if $r\equiv1$ mod 3 and $n-k=2\ell$, then for $\delta=0$, $1$, $2$, (\ref{recur})
becomes, respectively, $-\binom{\ell-1}m+\binom{\ell}m=\binom{\ell-1}{m-1}$, $\binom{\ell-1}m-\binom{\ell}m=-\binom{\ell-1}{m-1}$, and
$-\binom{\ell-1}m+0=-\binom{\ell-1}m$.

To initiate the induction we show that, mod 3,
\begin{equation}\label{T01}T_{0,1}(n,r)\equiv\begin{cases}2&n=2\cdot3^e\\
1&n=3^{e_1}+3^{e_2},\ 0\le e_1<e_2\\
r&n=3^e,\, e>0\\
r+1&n=1\\
0&\text{otherwise,}\end{cases}\end{equation}
and observe that the tabulated formulas for $k=0$ also equal (\ref{T01}).
The latter can be proved by considering separately $n=6t+d$ for $0\le d\le5$.
For example, if $d=3$, then $m=2t+1$, $\delta=0$, and $n-k=2(3t+1)+1$. For
$r\equiv0,1,2$, the tabulated value is, respectively, 0, $\binom{3t+1}{2t+1}$,
$-\binom{3t+1}{2t+1}$. Using Proposition \ref{Lagra}, one shows $\nu_3\bigl(\binom{3t+1}{2t+1}\bigr)=d_3(2t+1)-1$.
Thus the tabulated value in these cases is 0 mod 3 unless $2t+1$, hence $6t+3$, is a 3-power, and in this case
$\binom{3t+1}{2t+1}\equiv1$ mod 3.

To see (\ref{T01}), we note that
$$T_{0,1}(n,r)=\frac{3^{[(n-1)/2]}}{n!}S_1(n,r)$$
with $S_1$ as in Theorem \ref{mian}, and that, mod 3,
$$\frac{3^{[(n-1)/2]}}{n!}\equiv \begin{cases}1&n=3^{2e}\text{ or }3^e+3^{e+2k}\\
2&n=3^{2e+1},\,2\cdot3^e,\text{ or }3^e+3^{e+2k-1}\\
0&\text{otherwise.}\end{cases}$$
Thus, for example, mod 3, if $e>0$, then, using Theorem \ref{mian}
$$T_{0,1}(3^{2e},r)\equiv S_1(3^{2e},r)=(-1)^{(3^{2e}-1)/2}\bigl(\frac{3^{2e}+3}4+r\bigr)\equiv r,$$
in agreement with (\ref{T01}).
\end{pf}

\section{Proof of Theorem \ref{mainthm3}}\label{pf3sec}
In this section, we prove Theorem \ref{mainthm3}. We begin with a result, \ref{sum3}, which reduces much of the
analysis to evaluation of binomial coefficients mod 3.

\begin{defin}\label{taudef} For $\eps=\pm1$, let $\tau(n,k,\eps):=T_{k,1}(n,1)+\eps T_{k,1}(n,2)$, mod $3$.
\end{defin}
The following result is immediate from Theorem \ref{T1}.
\begin{prop} Let $n=3m+\delta$ with $0\le\delta\le2$. If $n-k=2\ell$, then, mod $3$, $\tau(n,k,-1)\equiv0$,
while $\tau(n,k,1)\equiv (-1)^\delta\binom{\ell-1}m$. If $n-k=2\ell+1$, then, mod $3$,
$$\tau(n,k,\eps)\equiv\begin{cases}0&\text{if }\delta=2\text{ or }\eps=1\text{ and }\delta=0\\
-\tbinom{\ell}m&\text{otherwise.}\end{cases}$$\label{tauval}
\end{prop}

The following result is a special case of Theorem \ref{phigen}, which is proved later.
\begin{thm}\label{sum3} Define
\begin{equation}\phi(n):=\sum\tbinom{n-1}k\tau([\tfrac n3],k,(-1)^{n-k-1})\in\Z/3.\label{fsum}\end{equation}
Then
$\nu_3(a_3(n-1,n))=s_3(n)$ if and only if $\phi(n)\ne0$.
\end{thm}

The following definition will be used throughout this section.
\begin{defin}\label{sparsedef} An integer $x$ is {\em sparse} if its base-3 expansion has no $2$'s or adjacent $1$'s.
The pair $(x,i)$ is {\em special} if $x$ is sparse and $i=x-\max\{3^{a_j}:\ 3^{a_j}\in x\}$.\end{defin}
\noindent Some special pairs are $(9,0)$, $(10,1)$, $(30,3)$, and $(91,10)$.

Lemma \ref{tech} will be used frequently. Its proof uses the following sublemma, which is easily proved.
\begin{sublem} Let $F_1(x,i)=(3x,3i)$ and $F_2(x,i)=(9x+1,9i+1)$. The special pairs are those that can be
obtained from $(1,0)$ by repeated application of $F_1$ and/or $F_2$.\end{sublem}
\noindent For example $(3^7+3^3+3,3^3+3)=F_1F_2F_2F_1F_1(1,0)$.

\begin{lem}\label{tech} Mod $3$,
\begin{enumerate}
\item If $x-i$ is even, then $\binom xi\binom{(3x-9i)/2}x\equiv0$;
\item If $x-i$ is odd, then $\binom xi\binom{(3x-9i-1)/2}x\equiv \begin{cases}1&\text{if $(x,i)$ special}\\
0&\text{otherwise;}\end{cases}$
\item If $x-i$ is odd, then $\binom xi\binom{(3x-9i-3)/2}x\equiv\begin{cases}1&\text{if $(x,i)$ special and }x\equiv0\, (3)\\
0&\text{otherwise.}\end{cases}$
\end{enumerate}
\end{lem}
\begin{pf} We make frequent use of (\ref{Lucas}).

(1) If $\binom xi\not\equiv0$, then $\nu_3(i)\ge\nu_3(x)$, but then the second factor is $\equiv0$ for a similar reason.

(2) Say $(x,i)$ satisfies $C$ if $\binom xi\binom{(3x-9i-1)/2}x\not\equiv0$. Note that $(1,0)$ satisfies $C$. We will show that
$(x,i)$ satisfies $C$ iff either $(x,i)=(3x',3i')$ and $(x',i')$ satisfies $C$ or $(x,i)=(9x''+1,9i''+1)$ and $(x'',i'')$ satisfies $C$.
The result then follows from the sublemma and the observation that the binomial coefficients maintain a value of 1 mod 3.

If $x=3x'$, then $\binom xi\not\equiv0$ implies $i=3i'$. Then
$$\tbinom{(3x-9i-1)/2}x\equiv\tbinom{(9x'-27i'-1)/2}{3x'}\equiv\tbinom{\frac12(9x'-27i'-3)+1}{3x'}\equiv\tbinom{(3x'-9i'-1)/2}{x'}.$$
If $x=3x'+1$, then $0\not\equiv\binom{\frac12(9x'-9i)+1}{3x'+1}$ implies $x'=3x''$. The product becomes
$\binom{9x''+1}i\binom{(3x''-i)/2}{x''}$. For this to be nonzero, $i$ cannot be $9i''$ by consideration of the second factor,
similarly to case (1). If $i=9i''+1$, the product becomes $\binom{x''}{i''}\binom{(3x''-9i''-1)/2}{x''}$, as claimed.
If $x=3x'+2$, a nonzero second factor would require the impossible condition $(9x'-9i+5)/2\equiv2$.

(3) To get nonzero, we must have $x=3x'$ then $i=3i'$. The product then becomes $\binom{x'}{i'}\binom{(3x'-9i'-1)/2}{x'}$,
which is analyzed using case (2).\end{pf}

Next we prove a theorem which, with \ref{sum3}, implies one part of the ``if" part of Theorem \ref{mainthm3}.
\begin{thm}\label{1thm} With $T$ as in Theorem \ref{mainthm3}, if
$n\in(3T+1)$
then $\phi(n)\ne0$.
\end{thm}
\begin{pf}  Define
$f_1(x)=\phi(3x+1)$.
The
lengthy proof breaks up into four cases, which are easily seen to imply the result, that
\begin{equation}\label{f1x}f_1(x)\ne0\text{ if }x\in T.\end{equation}
\begin{enumerate}
\item If $x$ is sparse, then $f_1(x)\ne0$.
\item For all $x$, $f_1(3x)=f_1(x)$.
\item If $x$ is not sparse and $x\not\equiv2$ mod 3, or if $x$ is sparse and $x\equiv1$ mod 3, then $f_1(3x+1)=\pm f_1(x)$.
\item If $x\equiv 0$ mod 3, then $f_1(3x+2)= f_1(x)$.
\end{enumerate}
Moreover, this inductive proof of (\ref{f1x}) will establish at each step that
\begin{eqnarray}\nonumber&&\text{if $\tbinom{3x}k\tau(x,k,(-1)^{x-k})\ne0$,
then $3x-k\equiv0$ (2)}\\  &&\text{unless $(3x,k)$ is special}.\label{stmt}\end{eqnarray}

\ss
{\bf Case 1}:  Let $x$ be sparse and $$3x=\sum_{j=1}^t 3^{a_j}$$
with  $a_j-a_{j-1}\ge2$ for  $2\le j\le t$.
Then
$$f_1(x)=\sum\tbinom{3x}{3i}\tau(x,3i,(-1)^{x-i}).$$
 We will show that
\begin{equation}\label{sps}\tbinom{3x}{3i}\tau(x,3i,(-1)^{x-i})=\begin{cases}-1&3i=3x-3^{a_t}\\
(-1)^j&3i=3x-3^{a_t}-3^{a_j},\,j\ge1\\
0&\text{otherwise}.\end{cases}\end{equation}
This will imply Case 1.

In the first case of (\ref{sps}), $(x,i)$ is special. If $x=3x'$, then $i=3i'$ with $(x',i')$ special, and we have
$$\tau(x,3i,-1)=-\tbinom{(3x'-9i'-1)/2}{x'}\equiv-1$$
by Lemma \ref{tech}.(2). If $x=3x'+1$, then $i=3i'+1$ with $(x',i')$ special.
Also, since $x$ is sparse, we must have $x'=3x''$ and then $i'=3i''$. Thus
$$\tau(x,3i,-1)=-\tbinom{(x-3i-1)/2}{x'}=-\tbinom{(3x''-9i''-1)/2}{x''}\equiv-1$$
by Lemma \ref{tech}.(2).

For the second case of (\ref{sps}), let $3i=3x-3^{a_t}-3^{a_j}$. This time $x-3i=2\ell$ with
$$\ell=\sum_{s=a_{t-1}}^{a_t-2}3^s+\cdots+\sum_{s=a_{j+1}}^{a_{j+2}-2}3^s+\sum_{s=a_{j-1}}^{a_j-2}3^s+\cdots
+\sum_{s=a_1}^{a_2-2}3^s+\sum_{s=a_1}^{a_{j+1}-2}3^s+2\cdot3^{a_1-1}.$$
Then $\ell-1$ is obtained from this by replacing $2\cdot3^{a_1-1}$ with $3^{a_1-1}+2\dstyle{\sum_{s= 0}^{a_1-2}3^s}$. Hence
$$\tau(x,3i,(-1)^{x-i})=(-1)^{\xbar}\tbinom{\ell-1}{[x/3]}\equiv 2^j\equiv(-1)^j.$$
Here we have used that for $\xbar=0,1$, we have $[\frac x3]=\dstyle{\sum_{j=\xbar+1}^t3^{a_j}}$.

We complete the argument for Case 1 by proving the third part of (\ref{sps}). The binomial coefficient $\binom{3x}{3i}$ is 0 unless
$3i=3x-3^{a_{j_1}}-\cdots-3^{a_{j_r}}$ with $j_1<\cdots<j_r$.  We must have $j_r=t$ or else $x-3i$ would be negative.
Hence $r>2$.
If $r=2w+1>1$ is odd, then
$$\tau(x,3i,(-1)^{x-i})=-\tbinom{\ell}{[x/3]}$$
with
$$2\ell+1=x-3i=\sum_{j\not\in\{j_1,\ldots,j_r\}}(3^{a_{j+1}-1}-3^{a_j})+\sum_{h=1}^w(3^{a_{j_{2h+1}}-1}+3^{a_{j_{2h}}-1})+3^{a_{j_1}-1},$$
and hence
$$\ell=\sum_{j\not\in\{j_1,\ldots,j_r\}}\sum_{i=a_j}^{a_{j+1}-2}3^i+\sum_{h=1}^w\biggl(3^{a_{j_{2h}-1}}+\sum_{i=a_{j_{2h}}-1}^{a_{j_{2h+1}}-2}
3^i\biggr)+\sum_{i=0}^{a_{j_1}-2}3^i.$$
Using (\ref{Lucas}), we see that $\binom{\ell}{[x/3]}\equiv0$ by consideration of position $a_{j_2}-2$.
A similar argument works when $r$ is even.
\ss

{\bf Case 2}: We are comparing
$$f_1(x)=\sum\tbinom{3x}{3i}\tau(x,3i,(-1)^{3x-3i})$$
with
$$f_1(3x)=\sum\tbinom{9x}{9i}\tau(3x,9i,(-1)^{9x-9i}),$$
mod 3. Clearly the binomial coefficients agree. Let $x=3y+\delta$ with $0\le\delta\le2$.

If $x-3i=2\ell$,
let $Q=(x-3i)/2$. We have
$$\tau(x,3i,1)=(-1)^\delta\tbinom{Q-1}y\equiv\tbinom{3Q-1}{3y+\delta}=\tau(3x,9i,1).$$
If $x-3i=2\ell+1$, let $Q=(x-3i-1)/2$. If $\delta\ne2$, we have
$$\tau(x,3i,-1)=-\tbinom Qy\equiv -\tbinom{3Q+1}{3y+\delta}=\tau(3x,9i,-1),$$
while if $\delta=2$, we have $\tau(x,3i,-1)=0$ by \ref{tauval}, and $\binom{3Q+1}{3y+\delta}=0$.

{\bf Case 3}:
Let $x=3y+\delta$ with $\delta\in\{0,1\}$.
Except for the single special term when $x$ is sparse,
we have $f_1(x)=\sum\binom{3x}{3i}\tau(x,3i,1)$, and will show
that
\begin{equation}f_1(3x+1)=\sum\tbinom{9x+3}{9i+3}\tau(3x+1,9i+3,1).\label{9x+3}\end{equation}
If $x-3 i=2\ell$, then $\tau(x,3i,1)=(-1)^\delta\binom{\ell-1}y$ and
$\tau(3x+1,9i+3,1)=-\binom{3\ell-2}{3y+\delta}\equiv-\tbinom{\ell-1}y$ since $\delta\ne2$. Thus $f_1(3x+1)=(-1)^{\delta+1}f_1(x)$.
To see that (\ref{9x+3}) contains all possible nonzero terms, note that terms $\binom{9x+3}{9i}\tau(3x+1,9i,(-1)^{x-i-1})$
contribute 0 to $f_1(3x+1)$ since the $\tau$-part is $-\binom{(3x-9i)/2}x\equiv0$ or $-\binom{(3x-9i-1)/2}x\equiv0$, since $(x,i)$ is not special.

If $x$ is sparse, the special term $(x,i)$ contributes $-1$ to $f_1(x)$.
If also $x\equiv1$ mod 3, then the corresponding term in (\ref{9x+3}) is $\tau(3x+1,9i+3,-1)$ with $x-i$ odd, equaling $-\binom{(3x-9i-3)/2}x\equiv-1$
by \ref{tech}.(3). That the terms added to each
are equal is consistent with $f_1(3x+1)=(-1)^{\delta+1}f_1(x)$.

{\bf Case 4}: Let $x=3y$. Ignoring temporarily the special term  when $x$ is sparse, we have
$f_1(x)=\sum\binom{3x}{3i}\tau(x,3i,1)$ and will show that $f_1(3x+2)=\sum\binom{9x+6}{9i+6}\tau(3x+2,9i+6,1)$. If $x-3i=2\ell$,
then
$$\tau(x,3i,1)\equiv\tbinom{\ell-1}y\equiv\tbinom{3\ell-3}{3y}\equiv\tau(3x+2,9i+6,1).$$
If the $9i+6$ in the sum for $f_1(3x+2)$ is replaced by $9i$ or $9i+3$, then the associated $\tau$ is 0,
for different reasons in the two cases.

We illustrate what happens to a special term $(x,i)$ when $x$ is sparse, using  the case $x=30$ and $i=3$.
It is perfectly typical. This term contributes $-1$ to $f_1(x)$. We will show that it also contributes $-1$ to $f_1(3x+2)$,
using $9i+3$ rather than $9i+6$, which is what contributed in all the other cases. The reader can check that for terms with
$k=9i+\la 0,3,6\ra$, the $\tau$-terms are, respectively
$$\tau(92,27,-1)=0,\,\tau(92,30,1)\equiv\tbinom{30}{30}\equiv1,\,\tau(92,33,-1)=0.$$
The binomial coefficient accompanying the case $i=30$ is $\binom{9\cdot30+6}{9\cdot3+3}\equiv2$.
\end{pf}

Next we prove a theorem, similar to \ref{1thm}, which, with \ref{sum3}, implies another part of the ``if" part of Theorem \ref{mainthm3}.
\begin{thm}\label{3thm} With $T$ as in Theorem \ref{mainthm3}, if
$n\in(9T+3)$
then $\phi(n)\ne0$.
\end{thm}
\begin{pf} We define $f_3(x)=\phi(9x+3)$   and write $2\in x$ to mean that a 2 occurs somewhere in the
3-ary expansion of $x$. We organize the proof into four cases, which imply the result.
\begin{enumerate}
\item If $2\not\in x$, then $f_3(x)\ne0$.
\item For all $x$, $f_3(3x)=f_3(x)$.
\item For all $x$, $f_3(9x+2)= f_3(x)$.
\item If $x$ is not sparse and $x\not\equiv2$ mod 3, then $f_3(3x+1)=(-1)^{\xbar+1}f_3(x)$.
\end{enumerate}
\ss
{\bf Case 1}: Let $9x=\dstyle{\sum_{i=1}^t 3^{a_i}}$ with $a_i>a_{i-1}$ and $a_1\ge2$. Let $i_0$ be the
largest $i\ge1$ such that $a_{i+1}-a_i=1$. Note that $x$ is sparse iff no such $i$ exists; let $i_0=1$
in this situation. For any $j$, let $p(j)$ denote the number of $i\le j$ for which $a_{i-1}<a_i-1$ or $i=1$.
We will sketch a proof that, mod 3,
\begin{equation}\label{pj}\tbinom{9x+2}k\tau(3x+1,k,(-1)^{x-k})\equiv\begin{cases}
1\cdot(-1)^{p(j)+1}&k=9x+2-3^{a_t}-3^{a_j},\,i_0\le j<t\\
2\cdot(-1)&k=9x+1-3^{a_t},\ n\text{ sparse}\\
0&\text{otherwise.}\end{cases}\end{equation}
We have written the values in a form which separates the binomial coefficient factor from the $\tau$ factor.
The binomial coefficient factor follows from (\ref{Lucas}). One readily verifies from (\ref{pj}) that the nonzero
terms in (\ref{fsum}) written in increasing $k$-order alternate between 1 and $-1$ until the last one which repeats
its predecessor. Thus the sum is nonzero.

The hard part in all of these is discovering the formula; then the verifications are straightforward,
and extremely similar to those of the preceding proof. We give one, that shows where $(-1)^{p(j)+1}$ comes from.

If $k=9x+2-3^{a_t}-3^{a_j}=2+3^{a_1}+\cdots+3^{a_{j-1}}+3^{a_{j+1}}+\cdots+3^{a_{t-1}}$, then $3x+1-k=2\ell+1$ with
$$\ell=\sum\Sb i=2\\ i\ne j+1\endSb^{t}\sum_{s=a_{i-1}}^{a_{i}-2}3^s+\sum_{s=0}^{a_{j+1}-2}3^s+\sum_{s=0}^{a_1-2}3^s.$$
We desire $\tau(3x+1,k,1)=-\binom\ell x$ with $x=\sum_{i=1}^t 3^{a_i-2}$. Note that $\ell$ has a $3^{a_i-2}$-summand
for each $i\ne j+1$ for which $a_{i-1}\ne a_i-1$, and another for each $i\le j+1$. Thus the 3-ary expansion of $\ell$ will have 0 in position
$a_i-2$, causing $\tau=0$, if $i>j+1$ and $a_i=a_{i-1}+1$. That explains the choice of $i_0$. If $j\ge i_0$, then $\binom\ell x$
from (\ref{Lucas}) has a factor $\binom 21$ in positions $i$ enumerated by $p(j)$.

\ss
{\bf Case 2}: If $x$ is sparse, the result follows from the proof of Case 1, and so we assume $x$ is not sparse. Then we are comparing
\begin{equation}\label{f3sum}f_3(x)=\sum\tbinom{9x+2}{9i+2}\tau(3x+1,9i+2,(-1)^{x-i}),\end{equation}
mod 3, with
\begin{equation}\label{f33sum}f_3(3x)=\sum\tbinom{27x+2}{27i+2}\tau(9x+1,27i+2,(-1)^{x-i}).\end{equation}
The binomial coefficients are clearly equal, mod 3. One can show that, for the other possible contributors to (\ref{f33sum}),
 $\tau(9x+1,27i+1,(-1)^{x-i+1})=0=\tau(9x+1,27i,(-1)^{x-i})$.
If $x-i$ is odd, the $\tau$-terms in (\ref{f3sum}) and (\ref{f33sum}) are 0, while if $x-i$ is even and $Q=\frac{x-3i}2$, then
$$\tau(3x+1,9i+2,1)\equiv-\tbinom{3Q-1}x\equiv-\tbinom{9Q-1}{3x}\equiv\tau(9x+1,27i+2,1).$$

{\bf Case 3}: If $x$ is not sparse,
we are comparing
$$f_3(x)=\sum\tbinom{9x+2}{9i+2}\tau(3x+1,9i+2,(-1)^{x-i})$$
with
\begin{equation}\label{3rsum}f_3(9x+2)=\sum\tbinom{9(9x+2)+2}{9(9i+2)+2}\tau(3(9x+2)+1,9(9i+2)+2,(-1)^{x-i}).\end{equation}
We will show below that no other terms can contribute to (\ref{3rsum}). Given this, then the binomial coefficients clearly agree, mod 3.

When $x-i$ is odd, the terms in both sums are 0, since they are of the form
$\tau(3m+1,3m+1-2\ell,-1)$.

Suppose $x-i$ is even. Let $Q=\frac{x-3i}2$. The first $\tau$ is $-\binom{3Q-1}x$, while the second is the negative of
$\tbinom{27Q-7}{9x+2}\equiv\tbinom{3Q-1}x$, as desired.

As a possible additional term in (\ref{3rsum}), if $k=9(9i+2)+2$ is replaced with $k=9(9i+\a)+\b$ with $0\le\a,\b\le2$,
which are the only ways to obtain a nonzero
binomial coefficient, then we show that the relevant $\tau$ is 0. Still assuming $x-i$ even, if $\a+\b$ is odd, then
we obtain $\tau(3m+1,3m+1-2\ell,-1)=0$, while if $\b=0$ and $\a\ne1$, then we obtain $\tau=\binom{3y}{9x+2}\equiv0$ for some $y$.
Finally, if $\b=2$ and $\a=0$, $$\tau=\tbinom{9(3x-9i)/2+2}{9x+2}\equiv\tbinom{(3x-9i)/2}{x}.$$
Since, in order to have $\tbinom{9(9x+2)+2}{9(9i+2)+2}\not\equiv0$, we must have $\nu_3(i)\ge\nu_3(x)$, we conclude
$\binom{(3x-9i)/2}x\equiv0$ mod 3. The case $x-i$ odd is handled similarly.

If $9x=3^{a_1}+\cdots+3^{a_t}$ is sparse and $9i=9x-3^{a_t}$, there is an additional term, $\binom{9x+2}{9i+1}\tau(3x+1,9i+1,1)\equiv1$,
in the sum for $f_3(x)$. The additional term in $f_3(9x+2)$ is
$$\tbinom{9(9x+2)+2}{9(9i+1)+2}\tau(3(9x+2)+1,9(9i+1)+2,1)\equiv\tbinom{\ell}m,$$
with $m=9x+2=2+3^{a_1}+\cdots+3^{a_t}$, and $2\ell+1=3(9x+2)+1-9(9i+1)-2$, so that
$$\ell=\sum_{j=2}^t\sum_{s=a_{j-1}+2}^{a_j}3^s+\sum_{s=2}^{a_1}3^s +2,$$
and so the additional term in $f(9x+2)$ is 1.

\ss
{\bf Case 4}: We first show
\begin{eqnarray*}(-1)^{\xbar+1}f_3(x)&=&(-1)^{\xbar+1}\sum\tbinom{9x+2}{9i+2}\tau(3x+1,9i+2,(-1)^{x-i})\\
&=&\sum\tbinom{27x+11}{27i+11}\tau(9x+4,27i+11,(-1)^{x-i})\\
&=&f_3(3x+1)\end{eqnarray*}
for $x\equiv0,1$ (3). Both $\tau$'s are 0 if $x-i$ is odd, while if $x-i$ is even and $x\equiv0,1$ (3), then
$$(-1)^{\xbar+1}\tau(3x+1,9i+2,1)\equiv(-1)^{\xbar}\tbinom{3Q+2}x\equiv-\tbinom{9Q+5}{3x+1}=\tau(9x+4,27i+11,1),$$
where $Q=(x-3i-2)/2$.

We must also show that $\binom{27x+11}k\tau(9x+4,k,(-1)^{x+1-k})\equiv0$ for $k\not\equiv11$ (27). When $k\equiv2$ (27),
the result follows from Lemma \ref{tech}. When $k\equiv0,9$ (27), $\tau$ is of the form $\binom{3A}{3x+1}\equiv0$.
\end{pf}

The ``if" part of Theorem \ref{mainthm3} when $n=3T'+2$ divides into two parts, Theorems \ref{5thm} and \ref{nthm},
noting that $3T'+2=(9T+2)\cup(9T'+5)$.
\begin{thm}\label{5thm} If $T$ is as in \ref{mainthm3} and $n\in(9T'+5)$, then $\phi(n)\ne0$.\end{thm}
\begin{pf}  Let $f_5(x)=\phi(9x+5)$.
We will prove that if $x\in T'$ then
\begin{equation}\label{351}f_5(x)=(-1)^{\xbar}f_3(x).\end{equation} With Theorem \ref{3thm}, this implies the result.

\ss
{\bf Case 1}: Assume $x$ not sparse and recall $x\not\equiv2$ (3). We show that, mod 3,
\begin{equation}\label{35}\tbinom{9x+4}k\tau(3x+1,k,(-1)^{9x+4-k})\equiv(-1)^{\xbar}\tbinom{9x+2}{k-2}\tau(3x+1,k-2,(-1)^{9x+4-k}).\end{equation}
Since $f_5(x)$ is the sum over $k$ of the LHS, and $(-1)^{\xbar}f_3(x)$ the sum over $k$ of the RHS, (\ref{351}) will follow when $x$ is not sparse.

We first deal with cases when the RHS of (\ref{35}) is nonzero. By the proof of \ref{3thm}, this can only happen
when $k-2=9i+2$, $\binom xi\not\equiv0$ mod 3,
and $x-i$ is even. Mod 3, we have $\binom{9x+4}{9i+4}\equiv\binom{9x+2}{9i+2}$ by (\ref{Lucas}). The two $\tau$'s in (\ref{35}) are,
with $Q:=\frac{3x-9i}2$, $-\binom{Q-2}x$ and $-\binom{Q-1}x$, respectively. Since $Q\equiv0$ mod 3,  these are
equal if $x\equiv0$ and negatives if $x\equiv1$.

We conclude the proof of (\ref{35}) by showing that other values of $k$ cause $\binom{9x+4}k\tau(3x+1,k,(-1)^{x-k})\equiv0$. If $k\not\equiv0,1,3,4$
mod 9, then $\binom{9x+4}k\equiv0$. If $k=9i+1$ or $9i+3$ and $x-i$ even, or if $k=9i$ or $9i+4$ and $x-i$ odd, then $\tau=0$ by \ref{tauval}.
If $k=9i$ and $x-i$ is even, then $\tau\equiv\binom{3x-9i}x\equiv0$. For $k=9i+1$ or $9i+3$ and $x-i$ odd, the result follows from Lemma \ref{tech}.

\ss
{\bf Case 2}: Assume $x$ is sparse. Let $9x=\dstyle{\sum_{j=1}^t 3^{a_j}}$ with $a_j-a_{j-1}\ge2$. We call $k=9i+d$, $d\in\{0,1,3,4\}$,
special if $(9x,9i)$ is special.
 The analysis of Case 1 shows that the $f_5$-sum over non-special values of $k$ equals $(-1)^{\xbar}$ times the $f_3$-sum over
non-special values of $k$.

We saw in (\ref{pj}) that the only special value of $k$ giving a nonzero summand for $f_3(x)$ is $k=9i+1$ (with $9i=9x-3^{a_t}$) and this summand is 1. We will show that if
$x\equiv1$ (3), then the only special value of $k$ giving a nonzero summand for $f_5(x)$ is $k=9i+1$, and it gives $-1$, while if $x\equiv0$ (3), both $k=9i+1$ and $k=9i+3$ give
summands of $-1$ for $f_5(x)$. This will imply the claim.

Recall $9i=9x-3^{a_t}$, and hence $x-i$ is odd. If $k=9i+\langle0,4\rangle$, then the $\tau$-factor is
$\tau(3x+1,9i+\langle0,4\rangle,-1)=0$. If $k=9i+\langle1,3\rangle$, the relevant term in $f_5(x)$ is
$$\tbinom{9x+4}{9i+\langle1,3\rangle}\tau(3x+1,9i+\langle1,3\rangle,1)=-\tbinom \ell x,$$
where $$\ell=\sum_{i=2}^t\sum_{s=a_{i-1}}^{a_i-2}3^s+\sum_{s=0}^{a_1-2}3^s+\langle0,-1\rangle.$$
Using (\ref{Lucas}), $\binom\ell x\equiv1$ in the $(9i+1)$-case, while in the $(9i+3)$-case
$$\binom\ell x\equiv\binom{(\sum_{s=0}^{a_1-2}3^s)-1}{3^{a_1-2}},$$
which is 0 if $x\equiv1$ (3), since then $a_1=2$, but is 1 if $x\equiv0$ (3) since then $a_1\ge3$.
\end{pf}

When $n\in(9T+2)$, the equality of $e_3(n-1,n)$ and $s_3(n)$ in Theorem \ref{mainthm3} comes not from $\nu_3(a_3(n-1,n))$, as it has in the other cases,
but rather from $\nu_3(a_3(n-1,n+1))$. To see this, we first extend Theorem \ref{sum3}, as follows.
\begin{thm}\label{phigen} If $N\ge n$, then $\nu_3(a_3(n-1,N))=s_3(n)$ iff $[N/9]=[n/9]$ and
$$\sum\tbinom{n-1}k\tau([\tfrac N3],k,(-1)^{n-k-1})\not\equiv0\,(3).$$
\end{thm}
\begin{pf} This is very similar to the proof, centered around (\ref{Stirl}), that Theorem \ref{helpthm} implies Theorem \ref{mainthm2}.
We have
$$0=(-1)^NS(n-1,N)N!=a_3(n-1,N)+3^{n-1}\sum(-1)^k\tbinom N{3k}k^{n-1}.$$
Thus $\nu_3(a_3(n-1,N))=s_3(n)$ iff $B\not\equiv0$ (3), where, mod 3,
\begin{eqnarray*}B&:=&{\tfrac1{[n/3]!}}\sum(-1)^k\tbinom N{3k}k^{n-1}\\
&\equiv&\sum_{d=1}^2{\tfrac1{[n/3]!}}\sum_{k\equiv d\,(3)}(-1)^k\tbinom N{3k}k^{n-1}\\
&\equiv&{\tfrac1{[n/3]!}}\sum_{d=1}^2(-1)^d\sum_j(-1)^j\tbinom N{9j+3d}\sum_{\ell}3^\ell j^\ell \tbinom{n-1}\ell d^{n-1-\ell}\\
&\equiv&{\tfrac1{[n/3]!}}\sum_{d=1}^2(-1)^d\sum_j(-1)^j\tbinom N{9j+3d}\sum_\ell 3^\ell\tbinom{n-1}\ell d^{n-1-\ell}\sum_i S(\ell,i)i!\tbinom ji\\
&\equiv&{\tfrac1{[n/3]!}}\sum_{d=1}^2(-1)^d\sum_j(-1)^j\tbinom N{9j+3d}\sum_i 3^i\tbinom{n-1}i d^{n-1-i}i!\tbinom ji\\
&\equiv&{\tfrac{[N/3]!}{[n/3]!}}\sum_i\tbinom{n-1}i(T_{i,2}(N,3)+(-1)^{n-1-i}T_{i,2}(N,6))\\
&\equiv&{\tfrac{[N/3]!}{[n/3]!}}\sum_i\tbinom{n-1}i(T_{i,1}([{\tfrac N3}],1)+(-1)^{n-1-i}T_{i,1}([{\tfrac N3}],2))\\
&=&{\tfrac{[N/3]!}{[n/3]!}}\sum_i\tbinom{n-1}i\tau([{\tfrac N3]},i,(-1)^{n-i-1}).\end{eqnarray*}
\end{pf}

The ``if" part of \ref{mainthm3} when $n\in (9T+2)$ now follows from Theorem \ref{phigen} and the following result.
\begin{thm}\label{nthm} If $T$ is as in \ref{mainthm3} and $n\in (9T+2)$, then
$$\sum\tbinom{n-1}k\tau([\tfrac{n+1}3],k,(-1)^{n-k-1})\not\equiv0\,(3).$$
\end{thm}
\begin{pf} We prove that for such $n$
\begin{equation}\label{23}\sum\tbinom{n-1}k\tau([{\tfrac{n+1}3}],k,(-1)^{n-k-1})\equiv
\sum\tbinom nk\tau([\tfrac{n+1}3],k,(-1)^{n-k})\end{equation}
and then apply Theorem \ref{3thm}. Note that the RHS is $\phi(n+1)$.

If $n=9x+2$ with $x$ not sparse, then the proof of \ref{3thm} shows that the nonzero terms of the RHS of
(\ref{23}) occur for $k=9i+2$ with $\binom xi\not\equiv0$ (3) and $x-i$ even. Now (\ref{23}) in this case follows from
\begin{equation}\label{x12}\tbinom{9x+1}{9i+1}\tau(3x+1,9i+1,1)\equiv-\tbinom xi\tbinom{(3x-9i-2)/2}x\equiv\tbinom{9x+2}{9i+2}\tau(3x+1,9i+2,1).
\end{equation}
One must also verify that no other values of $k$ contribute to the LHS of (\ref{23}); this is done by the usual methods.

If $n=9x+2$ with $x$ sparse, (\ref{x12}) holds unless
$(x,i)$ is special. For such $i$, the contribution to the RHS of (\ref{23}) using $k=9i+1$ is $2\cdot2\equiv1$.
The LHS of (\ref{23}) obtains contributions of $1\cdot2$ from both $k=9i$ and $k=9i+1$. Indeed both $\tau$'s equal
$-\binom{(3x-9i-1)/2}x\equiv-1$ by \ref{tech}.
\end{pf}
The ``if" part of Theorem \ref{mainthm3} is an immediate consequence of Theorems \ref{1thm}, \ref{3thm}, \ref{5thm}, and \ref{nthm}.
We complete the proof of Theorem \ref{mainthm3} by proving the following result.
\begin{prop} If $n$ is not one of the integers described in Theorem \ref{mainthm3}, then for all integers $N\ge n$ satisfying
$[N/9]=[n/9]$, we have
$$\sum\tbinom{n-1}k\tau([{\tfrac N3}],k,(-1)^{n-k-1})\equiv0\ (3).$$
\end{prop}
\begin{pf} We break into cases depending on $n$ mod 9, and argue by induction on $n$ with the integers ordered so that
$9x+3$ immediately  precedes $9x+2$.

\ss
{\bf Case 1}: $n\equiv0$ (9).
Let $n=9a$. If $[\frac N3]=3a$ or $3a+2$, then $\tau([\frac N3],k,(-1)^{a-1-k})=0$ by \ref{tauval}.
Now suppose $[\frac N3]=3a+1$. We show that for each nonzero term in
$$\sum_k\tbinom{9a-1}k\tau(3a+1,k,(-1)^{a-k-1})$$
with $a-k$ odd, the $(k+1)$-term is the negative of the $k$-term. Thus the sum is 0.

Both $\tau$'s equal $-\binom{(3a-k-1)/2}a$. Since $\binom{9a-1}k+\binom{9a-1}{k+1}
=\binom{9a}{k+1}$, the binomial coefficients are negatives of one another  unless $k+1=9t$ with $\binom at\not\equiv0$ (3).
Then $\nu(t)\ge\nu(a)$ and so $\binom{(3a-k-1)/2}a=\binom{(3a-9t)/2}a\equiv0$ (3), so the $\tau$'s were 0.

\ss
{\bf Case 2}: $n\equiv6,7,8$ (9). In these cases, $[N/9]=[n/9]$ implies $[N/3]=[n/3]$ and so we need not consider $N>n$.
By \ref{tauval}, $\tau(3x+2,k,(-1)^{x+1-k})=0$, which implies $\phi(9x+6)=0=\phi(9x+8)$. We have
$$\phi(9x+7)=\sum\tbinom{9x+6}k\tau(3x+2,k,(-1)^{x-k}).$$
This is 0 if $x-k$ is odd, while if $x-k$ is even, a summand is $\binom{9x+6}k\binom{(3x-k)/2}x$, which is 0 unless $k\equiv0$ (3) and hence $x\equiv0$ (3).
In the latter case, with $x=3x'$ and $f_1$ as in the proof of \ref{1thm}, we have
$\phi(n)=f_1(9x'+2)$, which, by Case 4 of the proof of \ref{1thm}, equals $f_1(3x')$, and this is 0 by induction unless $x'\in T$.

\ss
{\bf Case 3}: $n=9x+5$. If $x\equiv0,1$ (3), then $\phi(9x+5)=\pm\phi(9x+3)$ was proved in Case 1 of the proof of \ref{5thm}.
The induction hypothesis thus implies the result for $N=n$ in these cases. If $x=3y+2$, then
$$\phi(n)=\sum\tbinom{27y+22}k\tau(9y+7,k,(-1)^{y-k}).$$
The $k$-term is 0 if $y-k$ is odd, while if $y-k$ is even, $\tau=-\binom{(9y+6-k)/2}{3y+2}$. This is 0 unless $k\equiv2$ mod 3, but then $\binom{27y+22}k\equiv0$
(3). The $k$-term for $N=n+1$ is  nonzero iff the $k$-term in $\phi(n)$ is nonzero; this is true
because $\tau(3z+2,k,(-1)^{z-k})=\pm\tau(3z+1,k,(-1)^{z-k})$. Thus the sum for $N=n+1$ is 0 if $x\not\in T'$.

\ss
{\bf Case 4}: $n=9x+2$. Since, for $\eps=0$ or 2, $\tau(3x+\eps,k,(-1)^{x-k+1})=0$, we deduce that $\sum\binom{n-1}k\tau([\frac N3],k,(-1)^{n-k-1})=0$
for $N=n$ and $N=n+4$. For $N=n+1$, this is just the LHS of (\ref{23}). By (\ref{23}), it equals $\phi(n+1)$, which is 0 for $x\not\in T$ by the
induction hypothesis.

\ss
{\bf Case 5}: $n=9x+3$. Let $f_3(x)=\phi(9x+3)$. Let $x$ be minimal such that $x\not\in T$ and $f_3(x)$ has a nonzero summand.
By the proof of \ref{3thm}, $x$ is not 0 mod 3, 2 mod 9, 1 mod 9, or 4 mod 9.

If $x\equiv5$, 7, or 8 mod 9, then $f_3(x)$ has no nonzero summands. For example, if $x=9t+7$, the summands are
$\binom{81t+65}k\tau(27t+22,k,(-1)^{t-k-1})$. This is 0 if $t-k$ is even, while if $t-k$ is odd, the $\tau$-factor is
$\binom{(27t+21-k)/2}{9t+7}$. For this to be nonzero, we must have $k\equiv5$ or 7 mod 9, but these make the first factor 0.
Other cases are handled similarly.

One can show that for $\eps=0$, 1, 2,
$$\tau(3x+2,9i+\eps,(-1)^{x-i-\eps})=\pm\tau(3x+1,9i+\eps,(-1)^{x-i-\eps})\in\Z/3.$$
This implies that when we use $N=n+3$, nonzero terms will be obtained iff they were obtained for $n$.

\ss{\bf Case 6}: $n\equiv1,4$ mod 9. Let $f_1(x)=\phi(3x+1)$. By  the proof of Theorem \ref{1thm},
there can be no smallest $x\equiv0,1$ mod 3 which is not in $T$ and has $f_1(x)\ne0$. When using $N=n+2$ or, if $n\equiv1\ (9)$,
$N=n+5$, then the $k$-summands, $\binom{9x}k\tau(3x+1,k,(-1)^{x-k})$, $\binom{9x+3}k\tau(3x+2,k,(-1)^{x+1-k})$, and $\binom{9x}k\tau(3x+2,k,(-1)^{x-k})$,
are easily seen to be 0.
\end{pf}

\section{Discussion of Conjecture \ref{thm2}}\label{conjsec}
In this section we discuss the relationship between $\ebar_2(n)$,
$e_2(n-1,n)$, and $s_2(n)$. In particular, we discuss an approach to Conjecture \ref{thm2}, which
suggests that
the inequality $e_2(n-1,n)\ge s_2(n)$ fails by 1 to be sharp if $n=2^t$, while if $n=2^t+1$, it is sharp but the maximum value of $e_2(k,n)$
occurs for a value of $k\ne n-1$. The prime $p=2$ is implicit in this section; in particular, $\nu(-)=\nu_2(-)$ and $a(-,-)=a_2(-,-)$.

Although  our focus will be on the two families of $n$ with which Conjecture \ref{thm2} deals,
we are also interested, more generally, in the extent to which equality is obtained in each of the inequalities of
\begin{equation}s_2(n)\le e_2(n-1,n)\le \ebar_2(n).\label{inequs}\end{equation}
In Table \ref{bigtbl}, we list the three items related in (\ref{inequs}) for $2\le n\le38$, and also the smallest positive $k$ for
which $e_2(k,n)=\ebar_2(n)$.
 We denote this
as $k_{\text{max}}$, since it is the simplest $k$-value giving the maximum value of $e_2(k,n)$. Note that in this range $k_{\text{max}}$ always equals
$n-1$ plus possibly a number which is rather highly 2-divisible.

\begin{table}
\begin{center}
\caption{Comparison for (\ref{inequs}) when $p=2$}
\label{bigtbl}
\begin{tabular}{c|cccc}
$n$&$s_2(n)$&$e_2(n,n-1)$&$\ebar_2(n)$&$k_{\text{max}}$\\
\hline
2&1&1&1&1\\
3&2&2&2&2\\
4&4&4&4&3\\
5&5&5&6&$4+2^3$\\
6&6&6&8&$5+2^3$\\
7&7&8&8&6\\
8&10&11&11&7\\
9&11&11&12&$8+2^6$\\
10&12&12&14&$9+2^6$\\
11&13&13&15&$10+2^6$\\
12&15&15&15&11\\
13&16&18&18&12\\
14&17&21&21&13\\
15&18&22&22&14\\
16&22&23&23&15\\
17&23&23&24&$16+2^{11}$\\
18&24&24&26&$17+2^{11}$\\
19&25&25&28&$18+2^{11}$\\
20&27&27&28&$19+2^{11}$\\
21&28&28&28&20\\
22&29&29&30&$21+2^{10}$\\
23&30&31&31&22\\
24&33&34&34&23\\
25&34&36&38&$24+2^{16}$\\
26&35&37&40&$25+2^{16}5$\\
27&36&38&40&$26+2^{16}$\\
28&38&40&40&27\\
29&39&42&44&$28+2^{18}$\\
30&40&43&45&$29+2^{18}$\\
31&41&46&46&30\\
32&46&47&47&31\\
33&47&47&48&$32+2^{20}$\\
34&48&48&50&$33+2^{20}$\\
35&49&49&52&$34+2^{20}$\\
36&51&51&53&$35+2^{20}$\\
37&52&52&54&$36+2^{20}3$\\
38&53&53&56&$37+2^{20}7$
\end{tabular}
\end{center}
\end{table}

We return to more specific information leading to Conjecture \ref{thm2}.
To obtain the value of $\ebar_2(n)$, we focus on large values of $e_2(k,n)$.
For $n=2^t$ and $2^t+1$, this is done in the following conjecture, which implies Conjecture \ref{thm2}.
Note that $s_2(2^t)=2^t+2^{t-1}-2$, and $s_2(2^t+1)=2^t+2^{t-1}-1$. We employ the usual convention $\nu(0)=\infty$.
\begin{conj}\label{ethm} If $t\ge3$, then
$$e_2(k,2^t)\begin{cases}=\min(\nu(k+1-2^t)+2^t-t,\, 2^t+2^{t-1}-1)&\text{if }k\equiv-1\ (\mo\, 2^{t-1})\\
< 2^t+2^{t-1}-1&\text{if }k\not\equiv-1\ (\mo\, 2^{t-1});\end{cases}$$
$$e_2(k,2^t+1)\begin{cases}=\min(\nu(k-2^t-2^{2^{t-1}+t-1})+2^t-t,\, 2^t+2^{t-1})&\text{if }k\equiv0\ (\mo\, 2^{t-1})\\
<2^t+2^{t-1}&\text{if }k\not\equiv0\ (\mo\, 2^{t-1}).\end{cases}$$
\end{conj}
Note from this that conjecturally the smallest positive value of $k$ for which $e_2(k,n)$ achieves its maximum value is $n-1$ when
$n=2^t$ but is $n-1+2^{2^{t-1}+t-1}$ when $n=2^t+1$. The reason for this is explained in the next result,
involving a comparison of the smallest $\nu(a(k,j))$ values.

\begin{conj}\label{athm}
There exist odd 2-adic integers $u$, whose precise value varies from case to case,
such that
\begin{enumerate}
\item if $k\equiv-1\ (\mo\, 2^{t-1})$, then
\begin{eqnarray*}\nu(a(k,2^t+1))&=&\nu(k+1-2^t-2^{2^{t-1}+t-1}u)+2^t-t\\
\nu(a(k,2^t+2))&=&\nu(k+1-2^t-2^{2^{t-1}+t-2}u)+2^t-t+1\\
\nu(a(k,2^t+3))&=&\nu(k+1-2^t-2^{2^{t-1}+t-2}u)+2^t-t+1;
\end{eqnarray*}
\item if $k\equiv0\ (\mo\, 2^{t-1})$, then
\begin{eqnarray*}\nu(a(k,2^t+1))&=&\nu(k-2^t-2^{2^{t-1}+t-1}u)+2^t-t\\
\nu(a(k,2^t+2))&=&\nu(k-2^t-2^{2^{t-1}+t}u)+2^t-t+1\\
\nu(a(k,2^t+3))&=&\nu(k-2^t-2^{2^{t-1}+t-2}u)+2^t-t+2.
\end{eqnarray*}
\end{enumerate}
\noindent For other values of $j\ge 2^t$ $($resp. $2^t+1)$, $\nu(a(k,j))$ is at least as large as all the values appearing on the RHS above.
\end{conj}
\ni Note that, for fixed $j$, $\nu(a(k,j))$ is an unbounded function of $k$; it is the interplay among several values of $j$ which causes the
boundedness of $e_2(k,n)$ for fixed $n$.

We show now that Conjecture \ref{athm} implies the \lq\lq$=\min$"-part of Conjecture \ref{ethm}.
In part (1), the smallest $\nu(a(k,j))$ for $j\ge 2^t$ is
$$\begin{cases}\nu(k+1-2^t)+2^t-t&\text{if }\nu(k+1-2^t)\le2^{t-1}+t-2,\text{ using }j=2^t+1\\
2^t+2^{t-1}-1&\text{if }\nu(k+1-2^t)=2^{t-1}+t-1,\text{ using }j=2^t+2\\
2^t+2^{t-1}-1&\text{if }\nu(k+1-2^t)>2^{t-1}+t-1,\text{ using either.}\end{cases}$$
In part (2), the smallest $\nu(a(k,j))$ for $j\ge 2^t+1$ is
$$\begin{cases}\nu(k-2^t)+2^t-t&\text{if }\nu(k-2^t)\le2^{t-1}+t-2,\text{ using }j=2^t+1\\
2^t+2^{t-1}&\text{if }\nu(k-2^t)=2^{t-1}+t-1,\text{ using }j=2^t+2\\
2^t+2^{t-1}-1&\text{if }\nu(k-2^t)\ge2^{t-1}+t,\text{ using }j=2^t+1.\end{cases}$$

Conjecture \ref{athm} can be thought of as an application of Hensel's Lemma, following Clarke (\cite{Cl}).
We are finding the first few terms of the unique zero of the 2-adic function $f(x)=\nu(a(x,2^t+\eps))$ for
$x$ in a restricted congruence class.

\section{Relationships with algebraic topology}\label{topsec}
In this section, we sketch how the numbers studied in this paper are related to topics
in algebraic topology, namely James numbers and $v_1$-periodic homotopy groups.

Let $W_{n,k}$ denote the complex Stiefel manifold consisting of $k$-tuples of orthonormal vectors in $\C^n$,
and $W_{n,k}\to S^{2n-1}$ the map which selects the first vector. In work related to vector fields on spheres,
James (\cite{Ja}) defined $U(n,k)$ to be the order of the cokernel of
$$\pi_{2n-1}(W_{n,k})\to\pi_{2n-1}(S^{2n-1})\approx\Z,$$
now called James numbers. A bibliography of  many papers in algebraic topology devoted to studying these
numbers can be found in \cite{DSU}. It is proved in \cite{Lune} that
$$\nu_p(U(n,k))\ge\nu_p((n-1)!)-\et_p(n-1,n-k).$$
Our work implies the following sharp result for certain James numbers.
\begin{thm} If $p=2$ or $3$, $n$ is as in Theorems \ref{mainthm2} or \ref{mainthm3}, and $L$ is sufficiently
large, then
$$\nu_p(U((p-1)p^L+n,(p-1)p^L))=p^L-(p-1)[{\tfrac np}]-\nu_p(n)-\nbar.$$
\end{thm}
\begin{pf} We present the argument when $p=3$.
By \cite[4.3]{DSU} and \ref{mainthm3}, we have
$$\nu_3(U(2\cdot3^L+n,2\cdot3^L))=\nu_3((2\cdot3^L+n-1)!)-(n-1+\nu_3([n/3]!)).$$
Using Proposition \ref{Lagra}, this equals
$$\tfrac12(2\cdot3^L-n-1-d_3(n-1)-[\tfrac n3]+d_3([\tfrac n3])).$$
If $\nbar\ne0$ and $n=3m+\nbar$, this equals $3^L-2m-\nbar$, while if $n=3m$, we use $d_3(k-1)=d_3(k)-1+2\nu_3(k)$ to obtain
$3^L-2m-\nu_3(3m)$.
\end{pf}

The $p$-primary $v_1$-periodic homotopy groups of a topological space $X$, denoted $\vp_*(X)_{(p)}$ and defined in \cite{DM}, are a first approximation to
the $p$-primary actual homotopy groups $\pi_*(X)_{(p)}$. Each group $\vp_i(X)_{(p)}$ is a direct summand of some homotopy group $\pi_j(X)$.
It was proved in \cite{DSU} that for the special unitary group $SU(n)$, we have, if $p$ or $n$ is odd,
$$\vp_{2k}(SU(n))_{(p)}\approx\Z/p^{e_p(k,n)},$$
and $\vp_{2k-1}(SU(n))_{(p)}$ has the same order. The situation when $p=2$ and $n$ is even is slightly more complicated; it was discussed in
\cite{BDSU} and \cite{DP}. In this case, there is a summand $\Z/2^{e_2(k,n))}$ or $\Z/2^{e_2(k,n)-1}$ in $\vp_{2k}(SU(n))_{(2)}$.
From Theorems \ref{mainthm3} and \ref{mainthm2} we immediately obtain
\begin{cor} If $n$ is as in Theorem \ref{mainthm3} and $k\equiv n-1$ mod $2\cdot3^{s_3(n)}$, then
$$\vp_{2k}(SU(n))_{(3)}\approx\Z/3^{s_3(n)}.$$
If $n$ is as in Theorem \ref{mainthm2} and is odd, and $k\equiv n-1$ mod $2^{s_2(n)-1}$, then
$$\vp_{2k}(SU(n))_{(2)}\approx\Z/2^{s_2(n)}.$$
\end{cor}

We are especially interested in knowing the largest value of $e_p(k,n)$ as $k$ varies over all integers,
as this gives a lower bound for $\exp_p(SU(n))$, the largest $p$-exponent of any homotopy group of the space.
It was shown in \cite{DS1} that this is  $\ge s_p(n)$ if $p$ or $n$ is odd. Our work here immediately implies Corollary \ref{SUcor}
since $\vp_{2n-2}(SU(n))_{(p)}$
has $p$-exponent greater than $s_p(n)$ in these cases.
\begin{cor} \label{SUcor} If $p=3$ and $n$ is not as in \ref{mainthm3} or $p=2$ and $n$ is odd and not as in \ref{mainthm2}, then $\exp_p(SU(n))>s_p(n)$.
\end{cor}

  Table \ref{bigtbl} illustrates how we expect that $k=n-1$ will give almost
the largest group $\vp_{2k}(SU(n))_{(p)}$, but may miss by a small amount. There is much more that might be done along
these lines.

\def\line{\rule{.6in}{.6pt}}

\end{document}